\documentclass[12pt]{amsbook}

\usepackage{amsmath, amsthm, amssymb}
\input xy
\xyoption{all}
\usepackage{mathrsfs}
\usepackage{enumerate}
\usepackage{caption}
\usepackage{graphicx}
\usepackage{tikz}
\usepackage{enumitem}

\newtheorem{theorem}{Theorem}[chapter]
\newtheorem{proposition}[theorem]{Proposition}
\newtheorem{lemma}[theorem]{Lemma}

\theoremstyle{definition}
\newtheorem{definition}[theorem]{Definition}

\newtheorem{remark}[theorem]{Remark}
\newtheorem{example}[theorem]{Example}

\newcommand{\Z}{\ensuremath{\mathbb{Z}}}
\newcommand{\Q}{\ensuremath{\mathbb{Q}}}
\newcommand{\C}{\ensuremath{\mathbb{C}}}
\newcommand{\R}{\ensuremath{\mathbb{R}}}
\renewcommand{\P}{\ensuremath{\mathbb{P}}}

\newcommand{\Hom}{\ensuremath{\textrm{Hom}}}

\newcommand{\M}{\ensuremath{\overline{\mathcal{M}}}}
\renewcommand{\O}{\ensuremath{\mathcal{O}}}
\newcommand{\X}{\ensuremath{\mathcal{X}}}
\newcommand{\ev}{\ensuremath{\textrm{ev}}}

\renewcommand{\d}{\ensuremath{\partial}}

\newcommand{\D}{\ensuremath{\mathcal{D}}}

\newcommand{\Jac}{\ensuremath{\textrm{Jac}}}

\begin{document}

\pagestyle{empty}

\title{Mirror Symmetry Constructions}
\author{Emily Clader and Yongbin Ruan}

\maketitle
\thispagestyle{empty}

\setcounter{chapter}{-1}

\noindent The following notes are based on lectures by Yongbin Ruan during a special semester on the B-model at the University of Michigan in Winter 2014.  No claim to originality is made for anything contained within.

\vspace{2in}

\noindent The authors would like to thank Mark Shoemaker for his detailed and valuable comments.  Thanks are also due to Kentaro Hori for answering many of the authors' questions, and to all of the students in the course on which these notes are based, without whom many corrections and clarifications would not have been made.  The authors were partially supported by NSF RTG grant 1045119.

\vspace{2in}
\noindent
Please send questions or comments to cladere@ethz.ch.

\tableofcontents

\pagestyle{headings}
\setcounter{page}{1}
\pagenumbering{arabic}

\chapter*{Preface: The Idea of Mirror Symmetry}

The term ``mirror symmetry" is used to refer to a wide array of phenomena in mathematics and physics, and there is no consensus as to its precise definition.  In general, it refers to a correspondence that maps objects of a certain type--- manifolds, for example, or polynomials--- to objects of a possibly different type in such a way that the ``A-model" of the first object is exchanged with the ``B-model" of its image.  The phrases ``A-model" and ``B-model" originate in physics, and the various definitions of mirror symmetry arise from different ideas about the mathematical data that these physical notions are supposed to capture.

The Calabi-Yau A-model, for example, encodes deformations of the K\"ahler structure of a Calabi-Yau manifold, while the Calabi-Yau B-model encodes deformations of its complex structure.  There is also a Landau-Ginzburg A-model and B-model, which are associated to a polynomial rather than a manifold, and which are somewhat less geometric in nature.  The versions of mirror symmetry that we will consider in this course are:
\begin{itemize}[leftmargin=0.4in]
\item The Batyrev construction, which interchanges the Calabi-Yau A-model of a manifold and the Calabi-Yau B-model of its mirror manifold;
\item The Hori-Vafa construction, which interchanges the Calabi-Yau (or, more generally, semi-Fano) A-model of a manifold and the Landau-Ginzburg B-model of its mirror polynomial;
\item The Berglund-H\"ubsch-Krawitz construction, which interchanges the Landau-Ginzburg A-model of a polynomial and the Landau-Ginzburg B-model of its mirror polynomial.
\end{itemize}
In each case, mirror symmetry is a conjectural equivalence between the sets of data encoded by the two models.  In full generality it remains a conjecture, but many cases are known to hold.  The Calabi-Yau/Calabi-Yau mirror symmetry, for example, has been proven whenever the Calabi-Yau manifold $X$ is a complete intersection in a toric variety, and in some cases when $X$ is a complete intersection in a more general GIT quotient.

We should note that, in these notes, mirror symmetry will only be discussed as an interchange of cohomology groups (or ``state spaces") on the A-side and B-side.  At least in the Calabi-Yau case, however, both the A-model and the B-model are understood to capture much more data than these vector spaces alone.  The Calabi-Yau A-model, for example, can be encoded in terms of Gromov-Witten theory.  

The structure of the notes is as follows.  In Chapter \ref{toricgeometry}, we will review the fundamentals of toric geometry, which are necessary to explain the Batyrev construction.  Chapters \ref{BB}, \ref{HV}, and \ref{BHK} develop the three forms of mirror symmetry outlined above.  The Appendix reviews the basics of Chen-Ruan cohomology, a cohomology theory for orbifolds that is needed in order to define the state spaces of the Calabi-Yau A- and B-model, and that also provides a useful parallel to the definition of the states spaces in Landau-Ginzburg theory.

\chapter{Toric Geometry}
\label{toricgeometry}

Toric geometry is the study of a class of algebraic varieties whose structure is entirely encoded by combinatorial data.  Due to their simplicity, toric varieties provide a natural testing ground for many algebro-geometric ideas, and furthermore, they allow the statement of mirror symmetry to be expressed combinatorially.

The contents of this chapter are based heavily on Chapter 7 of \cite{Hori}.  Other good references for the basics of toric geometry include \cite{CLS} and \cite{Fulton}.

\section{Toric varieties and fans}

\begin{definition}
A {\bf toric variety} is a complex variety $X$ containing an algebraic torus $T := (\C^*)^r$ as an open dense subset, for which the action of $T$ on itself by multiplication extends to an action of $T$ on all of $X$.
\end{definition}

For example, complex projective space $\P^r$ is a toric variety.  The open dense torus is
\[\{[x_0: \cdots : x_r] \; | \; x_i \neq 0\} \subset \{[x_0: \cdots : x_r]\} = \P^r.\]

When the toric variety is normal (which will always be the case in what follows), it can be constructed from a combinatorial object known as a fan.

Let $N$ be a lattice\footnote{In the general theory of toric varieties, $N$ is allowed to be any abelian group of finite rank.  However, in these notes, we will assume that $N$ has no torsion, and hence is a lattice.}, a discrete subgroup of $\R^r$ for some $r$.  It follows that $N \cong \Z^r$, but by referring to $N$ abstractly as a lattice, we are {\it not} fixing an isomorphism.
Denote
\[N_{\R} := N \otimes_{\Z} \R \cong \R^r.\]

\begin{definition}
\label{cone}
A {\bf convex rational polyhedral cone} $\sigma \subset N_{\R}$ is a set of the form
\begin{equation}
\label{sigma}
\sigma = \{a_1v_1 + \cdots + a_k v_k \; | \; a_i \geq 0\}
\end{equation}
for $v_1, \ldots, v_k \in N$.  (We will sometimes say $\sigma$ is {\bf generated} by the vectors $\{v_1, \ldots, v_k\}$.)  A convex rational polyhedral cone is called {\bf strongly convex} if, furthermore, $\sigma \cap (-\sigma) = \{0\}$; that is, $\sigma$ does not contain any hyperplanes.

Since we will deal exclusively with strongly convex rational polyhedral cones in this course, they will be referred to simply as ``cones".
\end{definition}

See Figures \ref{fig1} and \ref{fig2} for examples.  From these illustrations, the notion of a ``face" of a cone should be intuitively clear.  To put it precisely, a {\bf face} of a cone $\sigma$ defined as in (\ref{sigma}) is a subset given by setting some collection of the $a_i$'s to zero.

\begin{figure}[t]
\begin{minipage}{.5\linewidth}
\centering
\begin{tikzpicture}
\draw[->] (-1.2,0) -- (2,0);
\draw[->] (0,-1.2) -- (0,2);
\draw[thick] (0,0) -- (2,2);
\draw[fill=gray, semitransparent] (0,0) -- (2,2) -- (0,2) -- (0,0);
\foreach \n in {-1,1}{
        \draw (\n,-3pt) -- (\n,3pt);
        \draw (-3pt,\n) -- (3pt,\n);
    }
\end{tikzpicture}
\caption{The cone generated by $(1,1)$ and $(0,1)$ in $N \cong \Z^2$.}
\label{fig1}
\end{minipage}
\begin{minipage}{.49\linewidth}
\centering
\begin{tikzpicture}
\draw[->] (-1.5,0) -- (1.5,0);
\draw[->] (0,-1.2) -- (0,2.5);
\draw[thick] (-1.5,0) -- (1.5, 0);
\draw[gray, fill=gray, semitransparent] (-1.5,0) -- (1.5,0) -- (1.5, 2.4) -- (-1.5, 2.4) -- (-1.5,0);
\foreach \n in {-1,1}{
        \draw (\n,-3pt) -- (\n,3pt);
        \draw (-3pt,\n) -- (3pt,\n);
    }
\draw (-3pt, 2) -- (3pt, 2);
\end{tikzpicture}
\caption{The cone generated by $(1,0)$ and $(-1,0)$ is not strongly convex.}
\label{fig2}
\end{minipage}
\end{figure}

\begin{definition}
A {\bf fan} is a collection $\Sigma$ of strongly convex rational polyhedral cones satisfying:
\begin{enumerate}
\item Each face of a cone in $\Sigma$ is also a cone in $\Sigma$.
\item The intersection of any two cones in $\Sigma$ is a face of each of them.
\end{enumerate}
\end{definition}

\begin{example}
\label{ex1}
Let $N \cong \Z^2$, and define the following vectors in $N_{\R} \cong \R^2$:
\begin{align*}
v_1&:=(1,0)\\
v_2&:= (0,1)\\
v_3&:=(-1,-1).
\end{align*}
Then there is a fan $\Sigma$ whose cones are generated by every proper subset of $\{v_1, v_2, v_3\}$, where by convention, the empty set of vectors generates the $0$-dimensional cone $\{0\}$.  This fan is illustrated in Figure \ref{P2}.
\end{example}

\begin{example}
\label{ex2}
Again, let $N \cong \Z^2$.  Then there is a fan $\Sigma_n$ whose cones are generated by proper subsets of $\{(1,0), (-1,n), (0,1), (0,-1)\}$.  This fan is pictured in Figure \ref{Hirz}.
\end{example}

\begin{figure}[t]
\begin{minipage}{.5\linewidth}
\centering
\begin{tikzpicture}
\draw[->] (-2,0) -- (2,0);
\draw[->] (0,-2) -- (0,2);
\draw[gray, fill=gray, semitransparent] (0,0) -- (2,0) -- (2,2) --(0,2) -- (0,0);
\draw[darkgray, fill=darkgray, semitransparent] (0,0) -- (0,2) -- (-2,2) --(-2,-2) -- (0,0);
\draw[lightgray, fill=lightgray, semitransparent] (0,0) -- (-2,-2) -- (2,-2) --(2,0) -- (0,0);
\draw[very thick] (0,0) -- (-2,-2);
\draw[very thick] (0,0) -- (0,2);
\draw[very thick] (0,0) -- (2,0);
\foreach \n in {-1,1}{
        \draw (\n,-3pt) -- (\n,3pt);
        \draw (-3pt,\n) -- (3pt,\n);
    }
\draw (1,3pt) node [above] {$v_2$};
\draw (-3pt,1) node [left] {$v_3$};
\draw (-0.9,-1) node [right] {$v_1$};
\end{tikzpicture}
\caption{The fan $\Sigma$ from Example \ref{ex1}.}
\label{P2}
\end{minipage}
\begin{minipage}{.49\linewidth}
\centering
\begin{tikzpicture}
\draw[->] (-2,0) -- (2,0);
\draw[->] (0,-2) -- (0,2);
\draw[lightgray, fill=lightgray, semitransparent] (0,0) -- (2,0) -- (2,2) --(0,2) -- (0,0);
\draw[gray, fill=gray, semitransparent] (0,0) -- (0,2) -- (-2,2) --(-2,-2) -- (0,0);
\draw[lightgray, fill=lightgray, semitransparent] (0,0) -- (-2,-2) -- (0,-2) -- (0,0);
\draw[gray, fill=gray, semitransparent] (0,0) -- (0,-2) -- (2,-2) -- (2,0) -- (0,0);
\draw[very thick] (0,0) -- (-2,-2);
\draw[very thick] (0,0) -- (0,2);
\draw[very thick] (0,0) -- (2,0);
\draw[very thick] (0,0) -- (0,-2);
\foreach \n in {-1,1}{
        \draw (\n,-3pt) -- (\n,3pt);
        \draw (-3pt,\n) -- (3pt,\n);
    }
\draw (1,3pt) node [above] {$v_1$};
\draw (3pt,1) node [right] {$v_3$};
\draw (-0.9,-1) node [right] {$v_2$};
\draw (3pt, -1) node [right] {$v_4$};
\end{tikzpicture}
\caption{The fan $\Sigma_n$ from Example \ref{ex2}.}
\label{Hirz}
\end{minipage}
\end{figure}

There are several constructions that yield a toric variety from the data of a fan.  Perhaps the most standard (see \cite{Fulton}) involves defining an affine variety for each cone in the fan and using the intersections of the cones to describe how to glue these affine varieties together.  This procedure is analogous to the way that one obtains projective space by gluing together the affine subsets on which a given coordinate is nonzero.  We will return to this perspective in Section \ref{charts}.

For the present, we will take a different approach to defining toric varieties from fans.  In the case of projective space, our approach will yield the quotient presentation
\[\P^r = (\C^{r+1} \setminus \{0\})/\C^*\]
as opposed to the decomposition into affines.

Fix a fan $\Sigma$, and let $\Sigma(1)$ denote the set of $1$-dimensional cones.  These are sometimes called ``edges" or ``rays"; explicitly, they are simply the cones generated by a single nonzero vector in $N$.

For each cone $\rho \in \Sigma(1)$, there is a primitive generator $v_{\rho} \in N$.  That is, $v_{\rho}$ generates $\rho$ in the sense of Definition \ref{cone}, and for all integers $k >1$ we have $\frac{1}{k} v_{\rho} \notin N$.  For convenience, choose an ordering $v_1, \ldots, v_n$ of these vectors, where $n = |\Sigma(1)|$.  By abuse of notation, we will often write $\Sigma(1) = \{v_1, \ldots, v_n\}$, identifying $1$-dimensional cones with their primitive generators.

Consider an affine space $\C^n$ with a coordinate $x_{\rho}$ for each $\rho \in \Sigma(1)$.  In accordance with the above ordering, we will sometimes write these coordinates as $x_1, \ldots, x_n$.  Inside this affine space, the {\bf discriminant locus} is defined as
\[Z(\Sigma) = \bigcup_{\substack{S \subset \Sigma(1)\\ S \text{ does not span a cone in }\Sigma}} V(I_S),\]
where $I_S$ is the ideal
\[I_S = (\{x_\rho \; | \; \rho \in S\}) \subset \C[x_1, \ldots, x_n].\]

Let
\[M = \Hom(N,\Z),\]
the dual lattice to $N$. Then there is a homomorphism
\[\phi: \Hom_{\Z}(\Z^{\Sigma(1)}, \C^*) \rightarrow \Hom_{\Z}(M, \C^*) = N\otimes_{\Z} \C^*\]
\[f \mapsto \bigg(m \mapsto \prod_{v \in \Sigma(1)} f(v)^{\langle m,v \rangle} \bigg),\]
where $\langle \cdot, \cdot \rangle$ denotes the natural pairing $M \otimes N \rightarrow \Z$.  Let $G = \ker(\phi)$.

There is an action of $G$ on $\C^n$ defined by
\[g(x_1, \ldots, x_n) = (g(v_1)x_1, \ldots, g(v_n)x_n)\]
for $g \in G$, where we identify $v_{\rho} \in \Sigma(1)$ with the corresponding standard basis vector for $\Z^{\Sigma(1)}$.  It is straightforward to check that this action preserves $Z(\Sigma)$.  Thus, one can define
\[X_{\Sigma} := (\C^n \setminus Z(\Sigma))/G.\]
This is the definition of the toric variety associated to a fan $\Sigma$.

Note that $X_{\Sigma}$ is indeed toric; the torus that acts is
\[T := (\C^*)^n/G \cong \Hom_{\Z}(M, \C^*) = N\otimes_{\Z} \C^* \cong (\C^*)^r,\]
acting by the quotient of the usual diagonal action of $(\C^*)^n$ on $\C^n$.  The first isomorphism in this chain is given by $\phi$.

For computations, it is helpful to describe $\phi$ concretely in coordinates.  Under the ordering of $\Sigma(1)$ as $\{v_1, \ldots, v_n\}$, one obtains coordinates $(t_1, \ldots, t_n)$ for $\Z^{\Sigma(1)} \cong \Z^n$.  Furthermore, choosing a basis $\{e^1, \ldots, e^r\} $ for $N$ with dual basis $\{e_1, \ldots, e_r\}$ gives an identification $\Hom(M, \C^*) \cong (\C^*)^r$.  In these coordinates,
\[\phi(t_1, \ldots, t_n) = \left(\prod_{i=1}^n t_i^{v_{i1}}, \ldots, \prod_{i=1}^n t_i^{v_{ir}}\right),\]
where
\[v_i = \sum_{i=1}^r v_{ij} e^j.\]
One consequence of this description is that $G$ can easily be computed by determining the linear relations satisfied by the vectors $v_1, \ldots, v_n$.  We will carry this out explicitly in Section \ref{charge} of this chapter.

Let us compute the toric varieties associated to the two fans described above.

\begin{example}
\label{projspace}
If $\Sigma$ is the fan from Example \ref{ex1}, then $r=2$ and $n=3$.  The discriminant locus is
\[Z(\Sigma) = V(x_1, x_2, x_3) = \{0\},\]
since the only subset of the vectors $v_1, v_2, v_3$ that does not span a cone is the entire set.  In the standard basis $\{(1,0), (0,1)\}$ for $N \cong \Z^2$, the homomorphism $\phi$ is
\[\phi: (\C^*)^3 \rightarrow (\C^*)^2\]
\[(t_1, t_2, t_3) \mapsto (t_1^{-1}t_2, t_1^{-1}t_3).\]
Thus,
\[G = \{(t,t,t) \; | \; t \in \C^*\} \cong \C^*.\]
We obtain
\[X_{\Sigma} = (\C^3 \setminus \{0\})/\C^* = \P^2.\]
To see how the torus $T = (\C^*)^2$ sits inside $\P^2$, reverse the homomorphism $\phi$ to write
\[(\C^*)^2 \cong (\C^*)^3/G\]
\[(\lambda_1,\lambda_2) \mapsto [1, \lambda_1, \lambda_2].\]
Thus, $T \subset \P^2$ as $\{[1:y:z] \; | \; y,z \neq 0\}$.
\end{example}

\begin{example}
Let $\Sigma_n$ be the fan from Example \ref{ex2}, in the special case where $r=2$ and $n=4$.  Then the discriminant locus is
\[Z(\Sigma_n) = V(x_1, x_2) \cup V(x_3, x_4),\]
and $\phi$ is given in the standard basis for $N = \Z^2$ by
\[(t_1, t_2, t_3, t_4) \mapsto (t_1t_2^{-1}, t_2^{-n}t_3t_4^{-1}).\]
The kernel is
\[G = \{(\lambda_1, \lambda_1, \lambda_1^n\lambda_2, \lambda_2) \} \cong (\C^*)^2.\]
To understand the variety
\[X_{\Sigma_n} = \bigg(\C^4 \setminus \big(\{x_1=x_2=0\} \cup \{x_3 = x_4=0\}\big)\bigg)/(\C^*)^2,\]
first take the quotient by the $\lambda_1$ factor inside $(\C^*)^2$.  This yields the quotient of $\C^4 \setminus \{0\}$ by $\C^*$ as
\[\frac{\C^4 \setminus \{0\}}{\lambda(x_1, x_2, x_3, x_4) \sim (\lambda x_1, \lambda x_2, \lambda^n x_3, x_4)},\]
which is the complement of the zero section in the total space of $\O_{\P^1}(n) \oplus \O_{\P^1}$.  Now, taking the quotient by the $\lambda_2$ factor in $(\C^*)^2$ projectivizes this bundle.  Thus,
\[X_{\Sigma_n} = \P(\O_{\P^1}(n) \oplus \O_{\P^1}),\]
which is the Hirzebruch surface $F_n$.  Via the isomorphism
\[T = (\C^*)^2 \cong (\C^*)^4/G\]
\[(\lambda_1, \lambda_2) \mapsto (\lambda_1, 1, \lambda_2, 1)\]
induced by $\phi$, one sees that $T$ sits inside $F_n$ as
\[\{(x,1,z,1) \; | \; x,z \neq 0\},\]
where the first two coordinates are the base coordinates and the second two are the fiber coordinates.
\end{example}

\begin{example}
\label{noncompact}
Consider a fan $\Sigma$ in $\Z^3$ with
\[\Sigma(1) = \{(1,0,1), (0,1,1), (-1,-1,1), (0,0,1)\}.\]
Specifically, each subset of $\Sigma(1)$ of size $1$ or $2$ generates a cone, and the $3$-dimensional cones are generated by all subsets of size $3$ except for $\{(-1,-1,1), (1,0,1), (0,1,1)\}$.

The homomorphism $\phi$ is easily computed in the standard basis:
\[\phi(t_1, t_2, t_3, t_4) = (t_1t_3^{-1}, t_2t_3^{-1}, t_1t_2 t_3 t_4).\]
Thus,
\[G = \{(t,t,t,t^{-3})\} \cong \C^*,\]
and one obtains
\[X_{\Sigma} = \O_{\P^2}(-3).\]
Note that this is a noncompact toric variety.  In combinatorial terms, the noncompactness of $X_{\Sigma}$ is reflected by the fact that the fan does not fill out the entire ambient vector space $\R^3$, leaving some directions ``free" from the constraints imposed by the $G$-action.
\end{example}

In general, given a fan $\Sigma$, the resulting toric variety $X_{\Sigma}$ is compact if and only if $\Sigma$ spans $\R^n$.  Such fans are referred to as {\bf complete}.

\section{The charge matrix}
\label{charge}

As mentioned previously, there is a straightforward way to read off the group $G$, and in particular the weights of the $G$-action on $\C^n$, from the equations satisfied by vectors in $\Sigma(1)$.

To explain this, we will need to assume that $G$ does not contain any finite groups as summands; thus, $G \cong (\C^*)^s$ for some $s$.  Then the embedding
\[(\C^*)^s \subset (\C^*)^n\]
induced by viewing $(\C^*)^s$ as the kernel of $\phi$ can be written in coordinates as
\begin{equation}
\label{Q}
(t_1, \ldots, t_s) \mapsto \left(\prod_{a=1}^s t_a^{Q_{1a}}, \ldots, \prod_{a=1}^s t_a^{Q_{na}}\right)
\end{equation}
for a matrix $Q$.

\begin{definition}
The matrix $Q_{ij}$ in equation (\ref{Q}) is called the {\bf charge matrix} for $X_{\Sigma}$.
\end{definition}

The terminology relates to a physical connection with the gauged linear sigma model, which we will discuss in Chapter \ref{HV}.  It should be noted that the representation of $Q$ as a matrix depends on an identification of $G$ with $(\C^*)^s$ and hence is not canonical.

By the definition of $G$ as the kernel of $\phi$, we have
\[\sum_{i=1}^n Q_{ia} v_{ik} = 0\]
for all $a= 1, \ldots, s$ and all $k =1, \ldots, r$.  In other words, the charge matrix gives $s$ linear relations
\[\sum_{i=1}^n Q_{ia} v_i = 0\]
satisfied by the vectors $v_1, \ldots, v_n$.

It follows that if $\Lambda$ is the lattice of linear relations on $\Sigma(1)$, then a representation of $G$ as $(\C^*)^s$ is equivalent to a choice of basis for $\Lambda$.  Having made such a choice, one has
\[X_{\Sigma} = (\C^n \setminus Z(\Sigma))/(\C^*)^s,\]
in which $(\C^*)^s$ acts by
\[(\lambda_1, \ldots, \lambda_s)\cdot (x_1, \ldots, x_n) = \left(\prod_{a=1}^s \lambda_a^{Q_{1a}} \cdot x_1, \ldots, \prod_{a=1}^s \lambda_a^{Q_{na}}\cdot  x_n\right).\]
Thus, searching for a basis of linear relations among the $1$-dimensional cones gives a quick way to read off the toric variety from its fan.

\begin{example}
Consider the fan $\Sigma$ from Example \ref{ex1}, for which $X_{\Sigma} = \P^2$.  Then the linear relations among the vectors $\{v_1, v_2, v_3\} = \{(1,0),(0,1), (-1,-1)\}$ are generated by the single relation
\[v_1 + v_2 + v_3 = 0.\]
Under the corresponding identification of $G$ with $\C^*$, then, we have
\[Q = \left( \begin{array}{c} 1 \\ 1 \\1 \end{array}\right),\]
and indeed,
\[\P^2 = (\C^3 \setminus \{0\})/\C^*,\]
in which $\C^*$ acts with weight $1$ in each factor.
\end{example}

\begin{example}
For the Hirzebruch surface $F_n$ described by the fan in Example \ref{ex2}, the linear relations among $v_1, \ldots, v_4$ are generated by
\[v_3 + v_4 = 0\]
and
\[v_1 + v_2 - n\cdot v_4 = 0.\]
This choice determines an identification of $G$ with $(\C^*)^2$ in which
\[Q = \left( \begin{array}{cc} 0 & 1 \\ 0 & 1 \\ 1 & 0 \\ 1 & -n\end{array}\right).\]
Note that, again, the columns give the weights of the two $\C^*$ actions in the quotient
\[X_{\Sigma_n} = \frac{\C^4 \setminus Z(\Sigma_n)}{(\C^*)^2}.\]
\end{example}

\section{Divisors on $X_{\Sigma}$}
\label{divisors}

Associated to each cone $\sigma \in \Sigma$, there is a subvariety
\[Z_{\sigma} = \{x_{\rho_1} = \cdots = x_{\rho_k} = 0\} \subset X_{\Sigma},\]
where $\rho_1, \ldots, \rho_k$ are the generators of $\sigma$.  It is easy to check that $Z_{\sigma}$ is $T$-invariant.  In fact, all of the $T$-invariant subvarieties of $X_{\Sigma}$ are of the above form.

The codimension of $Z_{\sigma}$ is equal to the dimension of $\sigma$, so elements $\sigma \in \Sigma(1)$ yield $T$-invariant {\it divisors}.  If $\sigma$ has primitive generator $\rho \in N$, we will write $D_{\rho}$ for the $T$-invariant divisor associated to $\sigma$.

The divisors $D_{\rho}$ in fact generate the group $A_{r-1}(X_{\Sigma})$ of Weil divisors modulo linear equivalence; in other words, any divisor on $X_{\Sigma}$ is linearly equivalent to a linear combination of $T$-invariant divisors.  Indeed, one can say even more: it is possible to use the combinatorics of the fan to determine when two linear combinations of the divisors $D_{\rho}$ are linearly equivalent.  The criterion involves the notion of the principal divisor associated to a character $m \in M$, which we explain below.

Recall that the torus $T$ is defined as $\Hom(M, \C^*)$.  Thus, an element $m \in M$ defines a holomorphic function on $T$ by evaluation.  Since $X_{\Sigma}$ is a compactification of $T$, this holomorphic function extends to a meromorphic function $f_m$ on all of $X_{\Sigma}$.  Let $(m)$ be the divisor of zeroes and poles of $f_m$.  One can check that $(m)$ is also given by the explicit combinatorial formula
\[(m) = \sum_{\rho \in \Sigma(1)} \langle m, v_{\rho} \rangle D_{\rho} \in \Z^{\Sigma(1)} \cong A_{r-1}(X_{\Sigma}).\]
We will not prove this formula, but let us check it in an example.

\begin{example}
Let $\Sigma$ again be the fan from Example \ref{ex1}, whose associated toric variety is $\P^2$.  Choose
\[m = (a,b) \in M \cong \Z^2.\]
A point $(\lambda_1, \lambda_2) \in T= (\C^*)^2$ can be viewed as the homomorphism
\[M \rightarrow \C^*\]
\[(p,q) \mapsto \lambda_1^p \lambda_2^q.\]
Thus, the function on $T$ associated to the point $(a,b)$ is
\[(\lambda_1, \lambda_2) \mapsto \lambda_1^a \lambda_2^b.\]
Recalling from Example \ref{projspace} that $T$ sits inside $\P^2$ as $\{[1:y:z]\}$, it is clear that the unique meromorphic extension $f_m$ of this function to all of $\P^2$ is
\[(x,y,z) \mapsto \left(\frac{y}{x}\right)^a \left(\frac{z}{x}\right)^b.\]
In particular, this confirms that $f_m$ has a zero of order $\langle m ,v_2 \rangle = a$ at $D_2 = \{y=0\}$ (or a pole of order $-a$, if $a$ is negative), and similarly for the other vectors in $\Sigma(1)$.
\end{example}

Using these principal divisors, we obtain an explicit description of $A_{r-1}(X_{\Sigma})$.

\begin{theorem}
Given Weil divisors $D = \sum a_{\rho} D_{\rho}$ and $D'= \sum a'_{\rho} D_{\rho}$ on $X_{\Sigma}$, the following are equivalent:
\begin{enumerate}
\item $D$ and $D'$ are linearly equivalent;
\item $D$ and $D'$ are homologically equivalent (that is, they define the same element of $H_{2(r-1)}(X_{\Sigma};\Z)$);
\item $D$ and $D'$ have the same associated line bundle;
\item $D$ and $D'$ differ by $(m)$ for some $m \in M$.
\end{enumerate}
\end{theorem}

The upshot of everything we have said in this section, then, is that we have a short exact sequence
\begin{equation}
\label{sequence}
0 \rightarrow M \rightarrow \Z^{\Sigma(1)} \rightarrow A_{r-1}(X_{\Sigma}) \rightarrow 0,
\end{equation}
where $A_{r-1}(X_{\Sigma}) \cong H_{2(r-1)}(X_{\Sigma}; \Z)$ is the Chow group of divisors modulo linear equivalence.  The first map in the sequence is
\[m \mapsto \bigg(\langle m, v_{\rho} \rangle\bigg)_{\rho \in \Sigma(1)},\]
while the second sends $(a_{\rho})_{\rho \in \Sigma(1)}$ to $\sum a_{\rho} D_{\rho}$.
This induces
\begin{equation}
\label{ses}
0 \rightarrow \Hom(A_{r-1}(X_{\Sigma}), \C^*) \rightarrow \Hom(\Z^{\Sigma(1)}, \C^*) \rightarrow \Hom(M, \C^*) \rightarrow 0,
\end{equation}
and the rightmost arrow in this sequence is precisely $\phi$.

There are a number of additional consequences of the exact sequences (\ref{sequence}) and (\ref{ses}) that we should mention.  First, there is now a canonical, basis-independent description of $G$ as
\[G \cong \Hom(A_{r-1}(X_{\Sigma}), \C^*).\]
In case $X_{\Sigma}$ is compact, this is equivalent via Poincar\'e duality to $G = H_2(X; \Z) \otimes_{\Z} \C^*$.

Second, the entries of the charge matrix can now be described in a more geometrically-motivated way.  Rather than expressing it in terms of linear relations on the fan $\Sigma$, the sequence (\ref{sequence}) shows that such relations can be viewed as elements of $A_{n-1}(X_{\Sigma})$.  Thus, a basis $C_1, \ldots, C_n$ for the Mori cone of effective curve classes yields a basis for the lattice $\Lambda$ of relations.  Associating the divisors $D_1, \ldots, D_n$ to the generators $v_1, \ldots, v_n$ of $\Sigma(1)$, one can show that
\[Q_{ia} = D_i \cdot C_a,\]
in which $\cdot$ denotes the intersection product.

\section{Characterizing toric orbifolds}
\label{charts}

An orbifold is a variety that is locally given as the quotient of an affine space by a finite group action.  In keeping with the idea that all of the structure of toric varieties is represented combinatorially, one can read off from the fan whether a particular toric variety is an orbifold.  To do so, we will need to briefly describe the construction of $X_{\Sigma}$ in terms of charts rather than quotients.

For each top-dimensional cone $\sigma$ in $\Sigma$, define
\[X_{\sigma} = \{x \in \C^n \setminus Z(\Sigma) \; | \; x_{\rho} \neq 0 \text{ for } \rho \notin \sigma\}/G \subset X_{\Sigma}.\]
In other words, $X_{\sigma} = X_{\Sigma_{\sigma}}$, where $\Sigma_{\sigma}$ is the fan consisting of $\sigma$ and its faces.  We claim that
\begin{equation}
\label{charteq}
X_{\Sigma} = \bigcup_{\sigma \text{ top-dimensional}} X_{\sigma}.
\end{equation}
To see this, choose a point $x \in X_{\Sigma}$, and let
\[S_x = \{\rho \; | \; x_{\rho} = 0\}.\]
Then $S_x$ must span a cone, for otherwise $x$ would lie in $Z(X_{\Sigma})$.  If $\sigma$ is the maximal cone containing the cone spanned by $S_x$, then $x \in X_{\sigma}$, which proves the claim.

Since $G$ is infinite, the decomposition (\ref{charteq}) does not immediately present $X_{\Sigma}$ as an orbifold, despite the fact that each of the local patches $X_{\sigma}$ is the quotient of an affine variety by a subgroup of $G$.  We require a criterion to ensure that the stabilizer of the $G$-action on each $X_{\sigma}$ is finite.

\begin{definition}
A fan is {\bf simplicial} if each of its top-dimensional cones gives a $\Q$-basis for $N$.
\end{definition}

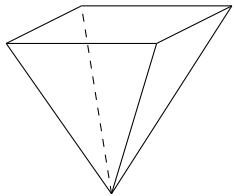
\begin{figure}[h]
\begin{tikzpicture}
\draw (0.4,0) -- (-1,2);
\draw (0.4,0) -- (1,2);
\draw (-1,2) -- (1,2);
\draw (-1,2) -- (0,2.5);
\draw (0,2.5) -- (2,2.5);
\draw (2,2.5) -- (1,2);
\draw (0.4,0) -- (2,2.5);
\draw[dashed] (0.4,0) -- (0,2.5);
\end{tikzpicture}
\caption{A non-simplicial cone.}
\label{conifold}
\end{figure}

\begin{example}
The cone in $\R^3$ depicted in Figure \ref{conifold} is not simplicial, since it has too many generators.  In fact, the toric variety corresponding to this fan is the conifold singularity $xy = uw$.
\end{example}

\begin{theorem}
A toric variety is an orbifold if and only if its fan is simplicial.  It is a (smooth) manifold if and only if each top-dimensional cone gives a $\Z$-basis for $N$.
\begin{proof}[Sketch of proof]

We will only give one direction of the proof, since we have explained how to construct a toric orbifold from a fan but not the reverse.  See Section 7.5 of \cite{Hori} for a discussion of the opposite implication.

Let $\Sigma$ be a simplicial fan, and let $\sigma$ be a top-dimensional cone.  Without loss of generality, we will write
\[\Sigma(1) = \{v_1, \ldots, v_n\}\]
and assume that $v_1, \ldots, v_r$ generate the $1$-dimensional faces of $\sigma$.

We claim that each $x \in X_{\sigma}$ is equivalent modulo the action of $G$ to a point with $x_{\rho} = 1$ for all $\rho \notin \sigma$.  Indeed, if $\rho \notin \sigma$, then $\rho \cup \sigma$ is linearly dependent, since $\sigma$ already has maximum dimension.  Thus, we have an equation
\[\sum_{i=1}^r a_i v_i  +a_{\rho}v_{\rho}= 0.\]
This implies that
\[(t^{a_1}, \ldots, t^{a_r},1, \ldots,  t^{a_{\rho}}, \ldots, 1) \in G\]
for any $t \in \C^*$, where the element $t^{a_{\rho}}$ is in the spot indexed by $v_{\rho}$.  Acting by this element, one can rescale $x_{\rho}$ to $1$ without changing any $x_{\rho'}$ for $\rho' \notin \sigma$ other than $\rho$.  Repeating this procedure for each $\rho \notin \sigma$ yields the claim.

It follows that
\begin{equation}
\label{stabilizer}
X_{\sigma} = \C^r/\text{stabilizer},
\end{equation}
where the stabilizer is the kernel of the restriction of $\phi$ to
\[\{x_{\rho} =1 \text{ for all } \rho \notin \sigma\} \subset (\C^*)^n.\]

Any element of $N$ yields an element of this stabilizer, since by assumption such an element may be written as
\[\ell_1 v_1 + \cdots + \ell_r v_r\]
with $\ell_i \in \Q$, and it is easily checked that
\[(e^{2\pi i \ell_1}, \ldots, e^{2\pi i \ell_r}, 1, \ldots, 1) \in \ker(\phi),\]
given that $\ell_1v_1 + \cdots + \ell_r v_r$ is integral.  This procedure produces the trivial element of $\ker(\phi)$ exactly when $\ell_i \in \Z$ for each $i$, so we find that the stabilizer in (\ref{stabilizer}) is isomorphic to $N/\Z\{v_1, \ldots, v_r\}$.  This is a finite group, and it is trivial exactly when $\{v_1, \ldots, v_r\}$ forms a $\Z$-basis for $N$.
\end{proof}
\end{theorem}

A particularly simple kind of orbifold is a {\bf global quotient}, which is the quotient of a smooth variety by a global finite group action.  In terms of fans, toric global quotients are described by passing to a finite-index sublattice.

Fix a fan $\Sigma$.  Suppose $N \subset N'$ is a sublattice of finite index such that the primitive\footnote{It should be noted that the definition of ``primitive" depends on the lattice.  In the lattice $\Z^2 \subset \R^2$, for example, the vector $(1,0)$ is a primitive generator for its ray.  The same ray, however, has primitive generator $(\frac{1}{2}, 0)$ if we instead consider the lattice $\frac{1}{2}\Z^2 \subset \R^2$.} generators of all of the top-dimensional cones of $\Sigma$ form integral bases for $N$.

Since $N_{\R} = N'_{\R}$, either of these lattices can be used with the same fan $\Sigma$ to define a toric variety.  The results, however, will be different.  Indeed, if
\[T = N \otimes_{\Z} \C^* \text{  and  } T' = N' \otimes_{\Z} \C^*,\]
then
\[X_{\Sigma,N} = \frac{X_{\Sigma,N}}{T/T} = \frac{X_{\Sigma, N}}{N/N}.\]
Though we will not prove this fact, let us see how it manifests in an example.

\begin{example}
\label{global}
Let $X_{\Sigma} = \P^2$, where the lattice is $N = \Z^2$.   Now, suppose we change the lattice to $N'= N + \Z\{(\frac{1}{3}, \frac{2}{3})\}$, which has the effect of adding two additional lattice points to each $1 \times 1$ square.  For example, the new lattice points in the square whose lower-left vertex is the origin are $(\frac{1}{3}, \frac{2}{3})$ and $(\frac{2}{3}, \frac{1}{3})$.

In the basis $\phi_1 = (\frac{2}{3}, \frac{1}{3})$, $\phi_2 = (\frac{1}{3}, \frac{2}{3})$ for $N'$, the generators of $\Sigma$ become
\[\{v_1, v_2, v_3\} = \{(2,-1), (-1,2), (-1,-1)\}.\]
As the generators of the one-dimensional cones no longer give integral bases for $N'$, it follows that $X_{\Sigma,N'}$ is an orbifold.  Indeed, one can check that $X_{\Sigma,N'} = [\P^2/\Z_3]$.
\end{example}

Conversely, given a global quotient $[X_{\Sigma}/H]$ of a toric variety, one can reconstruct the inclusion of lattices $N \subset N$ that defines it.  To do so, notice that $N$ can be read off from the torus $T = N \otimes_{\Z} \C^*$ as the lattice of $1$-parameter subgroups.  Although the torus $T \subset X_{\Sigma}$ descends to a torus $\overline{T} \subset [X_{\Sigma}/H]$, $1$-parameter subgroups of the latter need not lift to $1$-parameter subgroups of the former.  The new elements of $N'$ that were not present in $N$ are precisely the $1$-parameter subgroups of $\overline{T}$ that do not lift.  For instance, in Example \ref{global}, the subgroups
\[t \mapsto (1, t^{1/3}, t^{2/3})\]
and
\[t \mapsto (1, t^{2/3}, t^{4/3})\]
are well-defined in the quotient but {\it not} in the original variety $\P^2$, so these give the extra lattice points in $N'$.

\section{Toric resolutions}

Now that we have described toric orbifolds in some detail, let us discuss how their singularities can be resolved without leaving the toric setting.

\begin{definition}
A fan $\Sigma'$ is said to {\bf subdivide} $\Sigma$ if
\begin{enumerate}
\item $\Sigma(1) \subset \Sigma'(1)$;
\item Each cone of $\Sigma'$ is contained in some cone of $\Sigma$.
\end{enumerate}
\end{definition}

Suppose $\Sigma'$ subdivides $\Sigma$, and write
\[\Sigma(1) = \{\rho_1, \ldots, \rho_n\},\]
\[\Sigma'(1) = \{\rho_1, \ldots, \rho_n, \rho_{n+1}, \ldots, \rho_{m}\}.\]
Then the projection $\C^m \rightarrow \C^n$ determines a map $X_{\Sigma'} \rightarrow X_{\Sigma}$.  One can check that this map restricts to an isomorphism on the tori $T = N \otimes \C^*$ for $X_{\Sigma}$ and $X_{\Sigma'}$, so it is birational.  We refer to such a map as a {\bf toric resolution}.

\begin{figure}[ht]
\begin{tikzpicture}
\draw[->] (-1.2,0) -- (1.25,0);
\draw[->] (0,-1.2) -- (0,2.25);
\draw[gray, fill=gray, semitransparent] (0,0) -- (1.2,0) -- (1.2,2.2) -- (0,0);
\draw[very thick] (0,0) -- (1.2, 0);
\draw[very thick] (0,0) -- (1.2, 2.2);
\foreach \n in {-1,1}{
        \draw (\n,-3pt) -- (\n,3pt);
        \draw (-3pt,\n) -- (3pt,\n);
    }
\draw(-3pt, 2) -- (3pt, 2);
\end{tikzpicture}
\caption{The fan for $[\C^2/\Z_2]$.}
\label{blowup}
\end{figure}
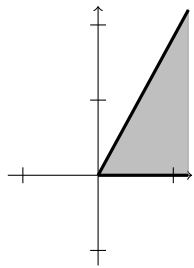

\begin{example}
Consider the toric orbifold $[\C^2/\Z_2]$.  By the procedure described in the previous section, this orbifold arises from the fan whose $1$-dimensional cones are
\[v_1 = (1,0) \in \Z^2,\]
\[v_2 = (0,1) \in \Z^2,\]
where the larger lattice is generated by $(1,0)$ and $(-\frac{1}{2}, \frac{1}{2})$.

In the basis $\{(1,0), (-\frac{1}{2}, \frac{1}{2})\}$, the coordinates of the vectors $v_i$ are $(1,0)$ and $(1,2)$, respectively.  Thus, the fan can equivalently be presented as in Figure \ref{blowup}.  If we add a one-dimensional cone generated by $v_3 = (1,1)$, subdividing the fan into two $2$-dimensional cones, then each of these cones has generators that give an integral basis for the lattice.  In other words, we have produced a smooth toric variety birational to $[\C^2/\Z_2]$. In fact, the resulting variety is $\O_{\P^1}(-2)$, which is the blowup of $[\C^2/\Z_2]$ at the origin.
\end{example}

\begin{figure}[t]
\begin{tikzpicture}
\draw[->] (-1.2,0) -- (1.25,0);
\draw[->] (0,-1.2) -- (0,4.25);
\draw[gray, fill=gray, semitransparent] (0,0) -- (1.2,0) -- (1.2,4.2) -- (0,0);
\draw[ultra thick] (0,0) -- (1.2, 0);
\draw[ultra thick] (0,0) -- (1.2, 4.2);
\draw[very thin] (0,0) -- (1.2, 3.2);
\draw[very thin] (0,0) -- (1.2, 2.2);
\draw[very thin] (0,0) -- (1.2, 1.2);
\foreach \n in {-1,1}{
        \draw (\n,-3pt) -- (\n,3pt);
        \draw (-3pt,\n) -- (3pt,\n);
    }
\draw(-3pt, 2) -- (3pt, 2);
\draw(-3pt, 3) -- (3pt, 3);
\draw(-3pt, 4) -- (3pt, 4);
\draw (-3pt,4) node [left] {$n+1$};
\end{tikzpicture}
\caption{The fan for $[\C^2/\Z_{n+1}]$ is depicted with bold lines.  The lighter lines show the rays that must be added to obtain the toric resolution.}
\label{n+1}
\end{figure}
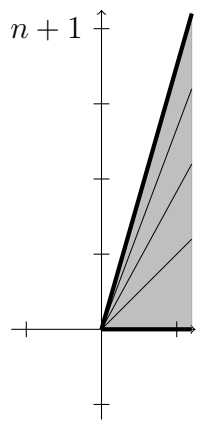

\begin{example}
Generalizing the above example, consider the global quotient $[\C^2/\Z_{n+1}]$, where the action is by
\[\omega\cdot (z_1, z_2) = (\omega z_1, \omega^n z_2)\]
for $\omega = e^{2\pi i/(n+1)}$.

Starting from the lattice $\Z^2$ that gives the toric variety $\C^2$, the extra lattice points that must be added to obtain $[\C^2/\Z_2]$ are generated by
\[ \left(\frac{n}{n+1}, \frac{1}{n+1}\right), \left(\frac{n-1}{n+1}, \frac{2}{n+1}\right), \ldots, \left(\frac{1}{n+1}, \frac{n}{n+1}\right).\]
Denote the corresponding vectors in $\R^2$ by $v_1, \ldots, v_n$; these generate the $1$-dimensional cones one must add to the fan for $[\C^2/\Z_2]$ in order to resolve its singularity.

In the refined lattice,
\[\left\{(1,0), \left(-\frac{1}{n+1},\frac{1}{n+1}\right)\right\}\]
is an integral basis for $\R^2$, and in this basis, $v_i$ has coordinates $(1,i)$.  Therefore, the fan for $[\C^2/\Z_{n+1}]$ and its blowup at the origin can be depicted as in Figure \ref{n+1}.
\end{example}

\begin{figure}[h]
\begin{tikzpicture}
\draw (0,0) -- (1,2);
\draw (0,0) -- (-1,2);
\draw (-1,2) -- (1,2);
\draw (-1,2) -- (0,2.25);
\draw (1,2) -- (0, 2.25);
\draw (-1,2) -- (0.2, 2.6);
\draw (1,2) -- (0.2,2.6);
\draw (0, 2.25) -- (0.2, 2.6);
\draw[dashed] (0,0) -- (0.2, 2.6);
\draw[dashed] (0,0) -- (0,2.25);
\end{tikzpicture}
\caption{The fan for $\O_{\P^2}(-3)$, which is a toric resolution of $[\C^3/\Z_3]$.}
\label{res2}
\end{figure}
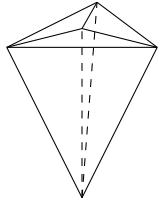

\begin{example}
Consider the global quotient $[\C^3/\Z_3]$, where $\Z_3$ acts with weight $1$ in each factor.  The fan for $\C^3$ is generated by the $1$-dimensional cones $(1,0,0), (0,1,0), (0,0,1)$ in the lattice $\Z^3$.  To obtain the quotient, one must add the lattice point $(\frac{1}{3}, \frac{1}{3}, \frac{1}{3})$ and its multiples.

A toric resolution of $[\C^3/\Z_3]$, then, is given by adding a $1$-dimensional cone generated by this new lattice point.  In the basis
\[\left\{\left(-\frac{1}{3}, \frac{2}{3}, -\frac{1}{3}\right), \left(-\frac{1}{3}, -\frac{1}{3}, \frac{2}{3}\right), \left(\frac{1}{3}, \frac{1}{3}, \frac{1}{3}\right)\right\}\]
for $\R^3$, which is integral under the refined lattice, the coordinates of the $1$-dimensional cones in the fan for $[\C^3/\Z_3]$ become
\[(-1,-1,1), (1,0,1), (0,1,1),\]
while the coordinates of the lattice point $(\frac{1}{3}, \frac{1}{3}, \frac{1}{3})$ become $(0,0,1)$.  Thus, the toric resolution can be depicted by the refinement of fans shown in Figure \ref{res2}.
\end{example}

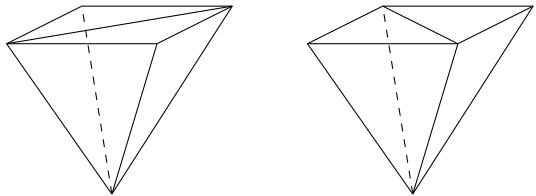
\begin{figure}[t]
\centering
\begin{tikzpicture}
\draw (0.4,0) -- (-1,2);
\draw (0.4,0) -- (1,2);
\draw (-1,2) -- (1,2);
\draw (-1,2) -- (0,2.5);
\draw (0,2.5) -- (2,2.5);
\draw (2,2.5) -- (1,2);
\draw (0.4,0) -- (2,2.5);
\draw[dashed] (0.4,0) -- (0,2.5);
\draw (-1,2) -- (2,2.5);

\draw (4.4,0) -- (3,2);
\draw (4.4,0) -- (5,2);
\draw (3,2) -- (5,2);
\draw (3,2) -- (4,2.5);
\draw (4,2.5) -- (6,2.5);
\draw (6,2.5) -- (5,2);
\draw (4.4,0) -- (6,2.5);
\draw[dashed] (4.4,0) -- (4,2.5);
\draw (4,2.5) -- (5,2);
\end{tikzpicture}
\caption{Two different toric resolutions of the conifold singularity.}
\label{resolutions}
\end{figure}

\begin{example}
Toric resolutions need not be unique.  For example, consider the (non-orbifold) fan consisting of the single cone shown in Figure \ref{conifold}, as well as its faces.  Specifically, the generators of the $1$-dimensional cones are $(1,0,0), (0,1,0), (0,0,1),$ and $(1,1,-1)$.  The resulting toric variety can be embedded as a hypersurface in $\C^4$ by the equation $\{xy = uv\}$, which is called the {\bf conifold singularity}.

There are two toric resolutions, given by adding the $2$-dimensional cones shown in Figure \ref{resolutions}.  The relationship between the two resulting toric varieties is known as a {\bf flop}.  It should be noted that neither is a blowup of the original conifold singularity, since they are given not by introducing new $1$-dimensional cones (which would correspond to adding divisors) but rather by adding a $2$-dimensional cone.
\end{example}

\section{Toric varieties from polytopes}

For projective toric varieties, an alternative combinatorial description can be given; instead of constructing the variety from the data of a fan, it is described via a polytope.

\begin{definition}
An {\bf integral polytope} in $M_{\R}$ is the convex hull of finitely many points in $M$.
\end{definition}

Given an integral polytope $\Delta \subset M_{\R}$ with integral points
\[\{m_0, \ldots, m_k\} = \Delta \cap M,\]
each $m_i$ can be viewed as a complex-valued function on the torus $T$, as explained in Section \ref{divisors}.  Let
\[f: T \rightarrow (\C^*)^{k+1} \subset \P^k\]
be the map defined by these $m_i$:
\[f(t) = [m_0(t): \cdots : m_k(t)].\]
It is easy to check that $f$ is an embedding, assuming that $\Delta$ is full-dimensional.

The toric variety associated to $\Delta$ is defined as
\[\P_{\Delta} := \overline{\text{Im}(f)}.\]
Since $f$ is injective, this contains $\text{Im}(f) \cong T$ as a dense open subset, and the fact that the $m_i$ are characters implies that the action of $T$ on itself by multiplication extends to all of $\P_{\Delta}$.  Hence, $\P_{\Delta}$ is indeed a toric variety.

Two observations can be made right away:
\begin{enumerate}
\item Shifting $\Delta$ by an integral point $m \in M$ has the effect of multiplying each coordinate of $f$ by the same complex number, which does not change the image in $\P^k$.  Thus, $\P_{\Delta}$ is independent of integral shifts, and in particular, we may assume without loss of generality that $0 \in M \cap \Delta$.
\item If the integral points of $\Delta$ satisfy an equation
\[\sum_{i=0}^k a_i m_i = 0\]
for integers $a_i$ with $\sum_{i=0}^k a_i = 0$, then the homogeneous coordinates $y_i$ of $\P_{\Delta}$ satisfy the equation
\[\sum_{a_i > 0} y_i^{a_i} = \sum_{a_i < 0} y_i^{-a_i}.\]
\end{enumerate}

We now have two procedures for obtaining normal toric varieties; the relationship between them is fairly straightforward.

\begin{definition}
Given an integral polytope $\Delta$, the {\bf normal fan} $\Sigma_{\Delta}$ has a cone $\sigma_F$ for each face $F$ of $\Delta$, defined by
\[\sigma_F = \{v \in N_{\R} \; | \; \langle m, v \rangle \leq \langle m', v \rangle \text{ for all } m \in F, m' \in \Delta\}.\]
\end{definition}

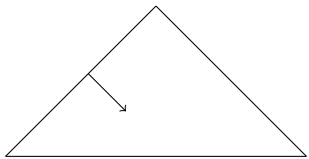
\begin{figure}[h]
\begin{tikzpicture}
\draw(-2,0) -- (2,0);
\draw (-2,0) -- (0,2);
\draw (2,0) -- (0,2);
\draw[->] (-0.9,1.1) -- (-0.4, 0.6);
\end{tikzpicture}
\caption{The inward-pointing normal vector to a face.}
\label{normalvector}
\end{figure}

In particular, the $1$-dimensional cones of the normal fan are generated by the integral normal vectors $v_F$ to the codimension-$1$ faces $F$ of $\Delta$, which are determined by the equation
\[F = \{m \; | \; \langle v_F, m \rangle = 0\}.\]
That is, if $N$ is identified with $M$ via a choice of basis, then $v_F$ is the inward-pointing integral normal vector to $F$--- see Figure \ref{normalvector}.

More generally, if $F$ is a codimension-$i$ face of $\Delta$, then $\sigma_F$ is an $i$-dimensional cone.  To determine $\sigma_F$, one writes $F$ as the intersection of a collection of codimension-$1$ faces; then $\sigma_F$ is generated by the integral normal vectors to the faces in this collection.

When $\P_{\Delta}$ is normal, there is an isomorphism
\[X_{\Sigma_{\Delta}} \cong \P_{\Delta}\]
given in homogeneous coordinates by
\begin{equation}
\label{fantopoly}
(x_1, \ldots, x_n) \mapsto \left(\prod_{i=1}^n x_i^{\langle m_0, v_i \rangle}, \ldots, \prod_{i=1}^n x_i^{ \langle m_k, v_i \rangle}\right).
\end{equation}
On the other hand, non-normal toric varieties can arise via polytopes, whereas the normal fan of a polytope is always defined.  (In combinatorial terms, this occurs when the lattice points $m_0, \ldots, m_k$ are a proper subset of $\Delta \cap M$ whose convex hull is nevertheless still $\Delta$.)  In this situation, the map (\ref{fantopoly}) is still an isomorphism on the torus, and hence is still birational, but it is not an isomorphism on the entire toric varieties; rather, $X_{\Sigma_{\Delta}}$ is a resolution of $\P_{\Delta}$.

\begin{figure}[t]
\begin{tikzpicture}
\draw[->] (-0.5,0) -- (2.5,0);
\draw[->] (0,-0.5) -- (0,2.5);
\draw[thick] (0,0) -- (2,0);
\draw[thick] (0,0) -- (0,2);
\draw[thick] (0,2) -- (2,0);
\draw[fill=gray, semitransparent] (0,0) -- (2,0) -- (0,2) -- (0,0);
\draw (2, -3pt) -- (2, 3pt);
\draw (-3pt, 2) -- (3pt, 2);
\end{tikzpicture}
\caption{The toric variety associated to this polytope is $\P^2$.}
\label{poly1}
\end{figure}
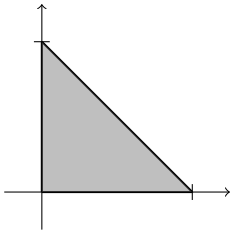

\begin{example}
Let $\Delta \subset \R^2$ be the polytope with vertices $(0,0), (1,0),$ and $(0,1)$, as shown in Figure \ref{poly1}.  The resulting map $f$ is
\[f: (\C^*)^2 \rightarrow \P^2\]
\[f(t_1, t_2) = [1: t_1 : t_2],\]
whose image is the standard embedding of the torus in $\P^2$.  Thus, $\P_{\Delta} = \P^2$, and after a shift of $\Delta$, the normal fan is the fan associated to $\P^2$.
\end{example}

\begin{example}
Scaling a polytope does not change the resulting toric variety (in particular, it does not affect the normal fan of the polytope), but it changes the embedding into projective space.  For example, if the polytope in the above example is scaled by a factor of two, then the map $f$ becomes
\[f: (\C^*)^2 \rightarrow \P^5\]
\[f(t_1, t_2) = [1:t_1:t_1^2: t_2: t_2^2: t_1t_2].\]
The closure of the image is a copy of $\P^2$ inside $\P^5$.
\end{example}

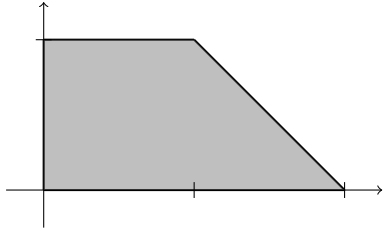
\begin{figure}[t]
\begin{tikzpicture}
\draw[->] (-0.5,0) -- (4.5,0);
\draw[->] (0,-0.5) -- (0,2.5);
\draw[thick] (0,0) -- (4,0);
\draw[thick] (4,0) -- (2,2);
\draw[thick] (2,2) -- (0,2);
\draw[thick] (0,2) -- (0,0);
\draw[fill=gray, semitransparent] (0,0) -- (4,0) -- (2,2) -- (0,2) -- (0,0);
\draw (2, -3pt) -- (2, 3pt);
\draw (4, -3pt) -- (4, 3pt);
\draw (-3pt, 2) -- (3pt, 2);
\end{tikzpicture}
\caption{The toric variety associated to this polytope is $F_2$.}
\label{poly2}
\end{figure}

\begin{example}
``Cutting off a corner" of a polytope adds a $1$-dimensional cone to the normal fan, so it yields a toric resolution.  For example, if the top corner of the polytope from the previous example is cut off, the resulting $\Delta$ is given in Figure \ref{poly2}.  The associated map is
\[f: (\C^*)^2 \rightarrow \P^4\]
\[f(t_1, t_2) = [1:t_2:t_1^2:t_2: t_2^2].\]
Either via the normal fan or via the polytope construction, it can be verified that the resulting toric variety is $F_2$.
\end{example}

When a toric variety $\P_{\Delta}$ is constructed out of a polytope, it is automatically equipped with an ample toric divisor $D$ defined as the pullback of the hyperplane class on projective space.  Thus, the polytope construction yields strictly more data than the fan method.  However, given a toric variety $X_{\Sigma}$ {\it together} with an ample toric divisor, it is possible to reconstruct the polytope $\Delta$ that yields this pair.\footnote{The case in which the divisor is ample but not very ample, and hence its associated morphism to projective space is not an embedding, is when the morphism $X_{\Sigma} \rightarrow \P_{\Delta}$ described in (\ref{fantopoly}) fails to be an isomorphism.}

Let $X_{\Sigma}$ be a toric variety equipped with an ample toric divisor $D$.  (The choice only matters up to linear equivalence, since a linearly equivalent divisor will yield a shift of the resulting polytope.)  Choosing $D$ amounts to specifying a morphism $X_{\Sigma} \rightarrow \P^k$ for which $\O(D)$ is the pullback of $\O(1)$, and hence the coordinate functions $x_0, \ldots, x_k$ on $\P^k$ yield sections of $\O(D)$.  Recall that there is an isomorphism
\[\left\{ \text{meromorphic functions } f \text{ on } X_{\Sigma} \; \big| \; D + (f) \geq 0\right\} \cong \Gamma(X_{\Sigma}, \O(D))\]
given by
\[f \mapsto f \cdot s_0,\]
in which $s_0$  is a global meromorphic function for which $(s_0 ) = D$.  Using this correspondence, the coordinate function $x_i$ yield meromorphic functions $f_i$ on $X_{\Sigma}$.  The restriction of each of these to the torus is a character of $T$, so it can be viewed as an element $m_i \in M$.  The polytope $\Delta$ is the convex hull of $m_0, \ldots, m_k$.

\begin{example}
\label{polyex}
Suppose we begin with the toric variety $\P^2$ and the toric divisor $D_0 = \{x_ 0 =0\}$.  Then $s_0 = x_0$ is a global meromorphic section whose divisor is $D_0$, so the functions $f_0, f_1, f_2$ corresponding to the coordinate sections of $D_0$ have
\[f_i \cdot x_0 = x_i,\]
and hence $f_i = x_i/x_0$.  In terms of the inhomogeneous coordinates $t_1 = \frac{x_1}{x_0}$ and $t_2= \frac{x_2}{x_0}$ on the torus, these are precisely
\[t_1^0t_2^0,\; t_1^1t_2^0, \; t_1^0t_2^1.\]
Thus, the polytope $\Delta$ is the convex hull of $(0,0)$, $(1,0)$, and $(0,1)$.
\end{example}

\begin{example}
Repeating the above example with the toric divisor $D_3$ instead of $D_1$ yields the convex hull of $(0,0)$, $(0,-1)$, and $(1,-1)$ as the polytope, which is indeed a shift of $\Delta$.
\end{example}

\begin{example}
\label{3dilate}
If, instead, the toric divisor is taken as $D_1 + D_2+D_3$, then the sections of $\O(3)$ are generated by homogeneous degree-three polynomials.  These can be dehomogenized by dividing each by $x_1x_2x_3$, at which point they can be expressed in terms of the torus coordinates $t_1 = \frac{x_2}{x_1}$ and $t_2 = \frac{x_3}{x_1}$ to read off the lattice points of the associated polytope.  For example, $x_1^3$ transforms to
\[\frac{x_1^3}{x_1x_2x_3} = \left(\frac{x_2}{x_1}\right)^{-1} \left(\frac{x_3}{x_2}\right)^{-1} = t_1^{-1}t_2^{-1},\]
so the point $(-1,-1)$ lies in the polytope.  Repeating this for each of the sections of $\O(3)$ reveals the polytope to be a threefold dilation of the polytope obtained in Example \ref{polyex}.
\end{example}

The above expression for the polytope is expressed more explicitly as follows:

\begin{proposition}
\label{polytope}
The integral points of the polytope associated to a toric variety $X_{\Sigma}$ with toric line bundle $\O\left(\sum_{\rho} a_{\rho} D_{\rho}\right)$ (in which $a_{\rho} \geq 0$) are
\begin{equation}
\label{arho}
\{m \in M \; | \; \langle m, v_{\rho} \rangle \geq -a_{\rho}\}.
\end{equation}
\begin{proof}
The previous definition shows that the integer points of the polytope consist of $m \in M$ for which $(m) + D \geq 0$.  This means that
\[\sum_{\rho} \langle m, v_{\rho} \rangle D_{\rho} + \sum_{\rho} a_{\rho} D_{\rho} \geq 0,\]
which is clearly equivalent to the description given in the statement of the proposition.
\end{proof}
\end{proposition}

In other words, the polytope $\Delta$ associated to $\big(X_{\Sigma}, \O(\sum_{\rho} a_{\rho} D_{\rho})\big)$ is bounded by the affine hyperplanes
\[F_{\rho} := \{m \in M_{\R} \; | \; \langle m, v_{\rho} \rangle = -a_{\rho}\}.\]
This fills out the correspondence between the fan and polytope perspectives on toric varieties.

\chapter{Batyrev Mirror Symmetry}
\label{BB}

The Batyrev mirror symmetry construction applies to toric varieties constructed out of a particular type of polytope, which we will discuss in the first subsection below.  Our presentation is based on Section 7.10 of \cite{Hori}; other key references for this material include \cite{Batyrev} and \cite{BB}.

\section{Reflexive polytopes}

\begin{definition}
A full-dimensional integral polytope $\Delta$ is {\bf reflexive} if there exist vectors $v_F \in N$ associated to each codimension-$1$ face $F$ of $\Delta$ such that
\[\Delta = \{m \in M_{\R} \; | \; \langle m, v_F \rangle \geq -1 \; \text{ for all } F\},\]
and if, furthermore, $0 \in \text{Int}(\Delta)$.
\end{definition}

A consequence of reflexivity is that $0$ is the only interior integral point of $\Delta$.

\begin{example}
\label{delta}
The polytope $\Delta$ appearing in Example \ref{3dilate} of the previous chapter, depicted in Figure \ref{3a} below, is reflexive.  The vectors $v_F$ associated to the edges are shown.
\end{example}

\begin{figure}[h]
\begin{tikzpicture}
\draw[->] (-1.5,0) -- (2.5,0);
\draw[->] (0,-1.5) -- (0,2.5);
\draw[thick] (-1,-1) -- (2,-1);
\draw[thick] (2,-1) -- (-1,2);
\draw[thick] (-1,2) -- (-1,-1);
\draw[fill=gray, semitransparent] (-1,-1) -- (2,-1) -- (-1,2) -- (-1,1);
\draw (2, -3pt) -- (2, 3pt);
\draw (1, -3pt) -- (1, 3pt);
\draw (-1, -3pt) -- (-1, 3pt);
\draw (-3pt, 2) -- (3pt, 2);
\draw (-3pt, 1) -- (3pt, 1);
\draw (-3pt, -1) -- (3pt, -1);
\draw (1.25,-1) node [below] {$v_{F_1} = (0,1)$};
\draw (0.5,0.75) node [right] {$v_{F_2} = (-1,-1)$};
\draw (-1,0.6) node [left] {$v_{F_3}= (1,0)$};
\end{tikzpicture}
\caption{A reflexive polytope.}
\label{3a}
\end{figure}
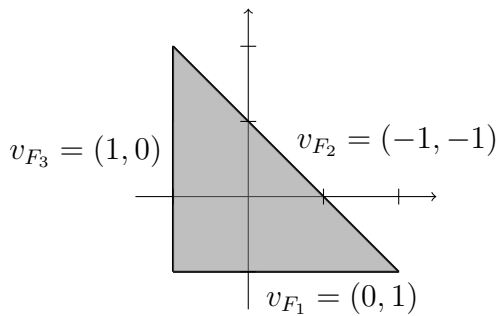

Proposition \ref{polytope} above shows that in the reflexive case, the divisor in $X_{\Sigma_{\Delta}}$ determined by the polytope is
\[D = \sum_{\rho \in \Sigma_{\Delta}(1)} D_{\rho}.\]
More explicitly, one can check (Theorem 8.2.3 of \cite{CLS}) that $D$ is the anticanonical divisor of $\X_{\Sigma_{\Delta}}$.  By exploiting this connection to the canonical divisor, one can prove the geometric meaning of reflexivity.

\begin{theorem}
\label{reflexive}
A full-dimensional lattice polytope $\Delta$ is reflexive if and only if $\P_{\Delta}$ is Gorenstein and Fano.
\end{theorem}

Recall that the Gorenstein condition implies that the canonical bundle extends across the singularities of $\P_{\Delta}$, and hence forms an honest line bundle on the entire variety.  Given this, the Fano condition says that the resulting bundle is ample.

Rather than prove the theorem, we trace its manifestation in the particular case of weighted projective space.

\begin{example}
By definition, weighted projective space is the quotient
\[\P(c_0, \ldots, c_n) = \frac{(\C^{n+1} \setminus \{0\})}{\C^*},\]
where $\C^*$ acts by
\[\lambda(z_0, \ldots, z_n) := (\lambda^{c_0}z_0, \ldots, \lambda^{c_n}z_n).\]

Weighted projective space is easily seen to be a toric variety.  In order to construct its fan, one must find $v_0, \ldots, v_n$ satisfying a relation
\[c_0 v_0 + \cdots + c_n v_n=0.\]

Assume, for simplicity, that $c_0=1$.  Then the vectors $v_i$ can be taken to be the following:
\[\left(\begin{array}{c}v_1 \\ \vdots \\ v_n \\ v_0 \end{array}\right) = \left(\begin{array}{cccc} 1 & 0 & \cdots & 0\\ 0 & 1 & \cdots & 0\\  \vdots & \vdots &  & \vdots\\ 0 & 0 & \cdots & 1\\ -c_1 & -c_2 & \cdots & -c_n\end{array}\right).\]
For example, the fan for $\P(1,c_1, c_2)$ is shown in Figure \ref{weightedproj}.

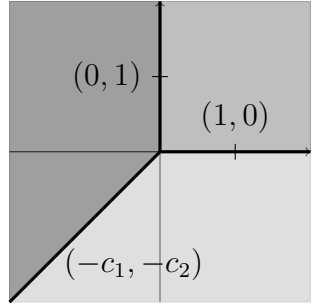
\begin{figure}[t]
\begin{tikzpicture}
\draw[->] (-2,0) -- (2,0);
\draw[->] (0,-2) -- (0,2);
\draw[gray, fill=gray, semitransparent] (0,0) -- (2,0) -- (2,2) --(0,2) -- (0,0);
\draw[darkgray, fill=darkgray, semitransparent] (0,0) -- (0,2) -- (-2,2) --(-2,-2) -- (0,0);
\draw[lightgray, fill=lightgray, semitransparent] (0,0) -- (-2,-2) -- (2,-2) --(2,0) -- (0,0);
\draw[very thick] (0,0) -- (-2,-2);
\draw[very thick] (0,0) -- (0,2);
\draw[very thick] (0,0) -- (2,0);
\draw(-3pt, 1) -- (3pt, 1);
\draw (1,-3pt) -- (1,3pt);
\draw (1,3pt) node [above] {$(1,0)$};
\draw (-3pt,1) node [left] {$(0,1)$};
\draw (-1.42,-1.5) node [right] {$(-c_1, -c_2)$};
\end{tikzpicture}
\caption{The fan for $\P(1, c_1, c_2)$.}
\label{weightedproj}
\end{figure}

Let $\Delta$ be the polytope $\Delta$ associated to $\P(1,c_1, \ldots, c_n)$ and its anticanonical divisor.  If $\Delta$ is reflexive, then it is defined by the inequalities $\langle m, v_i \rangle \geq -1$  for $0 \leq i \leq n$.  In particular, vertices of $\Delta$ would correspond to points where all but one of these is an equality:
\[\langle m_, v_j\rangle = -1 \text{ for } j \neq i.\]
These equations for $i \neq 0$ imply that any point $m = (m_0, \ldots, m_n) \in \Delta$ satisfies
\[m_j = -1 \text{ for } j \neq i,\]
\[-c_1m_1 - \cdots -c_n m_n = -1,\]
which implies that
\[c_ix_i = 1 + \sum_{j \neq 0,i} c_j = \sum_{j \neq i} c_j.\]
This equation must have an integral solution in order for such a polytope $\Delta$ to exist; thus, one must have
\[c_i \left| \sum_{j=0}^n c_j\right.\]
\end{example}

We have only proved one implication of Theorem \ref{reflexive}, and we have made the simplifying assumption that $c_0=1$, but in fact, the analysis holds more generally:

\begin{proposition}
The following are equivalent:
\begin{enumerate}
\item $\P(c_0, \ldots, c_n)$ is Gorenstein.
\item The polytope $\Delta$ associated to $\P(c_0, \ldots, c_n)$ with its anticanonical bundle is reflexive.
\item $c_i \bigg| \sum_{j=0}^n c_j$ for each $i$.
\end{enumerate}
Moreover, when these are satisfied, the polynomial
\[x_0^{d/c_0} + \cdots + x_n^{d/c_n}\]
with $d:=\sum_{j=0}^n c_j$ defines a Calabi-Yau hypersurface of $\P(c_0, \ldots, c_n)$.
\end{proposition}

The last statement follows from the adjunction formula.  Polynomials of the form $x_0^{a_0} + \cdots + x_n^{a_n}$ are referred to as {\bf Fermat}.  Thus, the result implies that a weighted projective space is Gorenstein if and only if it contains a Calabi-Yau hypersurface defined by a Fermat polynomial.

\section{Polar duality}

\begin{definition}
Suppose that $\Delta \subset M_{\R}$ is a full-dimensional lattice polytope containing $0$ as an interior point.  Then the {\bf polar dual} of $\Delta$ is
\[\Delta^{\circ} := \{ n \in N_{\R} \; | \; \langle m, n \rangle \geq -1 \text{ for all } m \in \Delta\}.\]
\end{definition}

It is straightforward to check that $(\Delta^{\circ})^{\circ} = \Delta$, justifying the name ``duality", and furthermore, that $\Delta$ is reflexive if and only if $\Delta^{\circ}$ is reflexive.

\begin{figure}[t]
\begin{tikzpicture}
\draw[->] (-1.5,0) -- (2.5,0);
\draw[->] (0,-1.5) -- (0,2.5);
\draw[thick] (-1,-1) -- (1,0);
\draw[thick] (1,0) -- (0,1);
\draw[thick] (0,1) -- (-1,-1);
\draw[fill=gray, semitransparent] (-1,-1) -- (1,0) -- (0,1) -- (-1,-1);
\draw (2, -3pt) -- (2, 3pt);
\draw (1, -3pt) -- (1, 3pt);
\draw (-1, -3pt) -- (-1, 3pt);
\draw (-3pt, 2) -- (3pt, 2);
\draw (-3pt, 1) -- (3pt, 1);
\draw (-3pt, -1) -- (3pt, -1);
\end{tikzpicture}
\caption{The polar dual of the polytope in Figure \ref{3a}.}
\label{polar}
\end{figure}
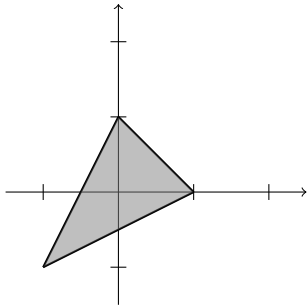

\begin{example}
The polar dual of the polytope $\Delta$ in Figure \ref{3a} is shown in Figure \ref{polar}.  Notice that the vectors $v_F$ for $\Delta$ become the vertices of $\Delta^{\circ}$.
\end{example}

One can see from this example that the normal fan $\Sigma_{\Delta}$ of a polytope $\Delta$ can also be obtained by taking cones over faces of $\Delta^{\circ}$.  This is indeed a general feature of polar duality, and will be important when understanding the relationship between $\P_{\Delta}$ and $\P_{\Delta^{\circ}}$ encoded by Batyrev mirror symmetry.

\section{Batyrev Mirror Symmetry}

If $\P_{\Delta}$ is the toric variety associated to a reflexive polytope $\Delta$, then any hypersurface defined by the vanishing of a generic section of the anti-canonical bundle $\O\left(\sum_{\rho} \D_{\rho}\right)$ will automatically be Calabi-Yau by the adjunction formula.  Varying the section yields a family of Calabi-Yau hypersurfaces, denoted by $X \subset \P_{\Delta^{\circ}}$.

\begin{definition}
Let $\Delta$ be a reflexive polytope  The {\bf Batyrev mirror} of the family of Calabi-Yau hypersurfaces $X \subset \P_{\Delta}$ as above is the family $X^{\circ} \subset \P_{\Delta^{\circ}}$ of hypersurfaces defined by the vanishing of a generic section of the anticanonical bundle of $\P_{\Delta^{\circ}}$.
\end{definition}

In what sense are these two families of Calabi-Yau hypersurfaces mirror to one another?  One fairly straightforward yet important fact is the {\bf monomial-divisor correspondence}, which is a bijection
\[\left\{ \begin{array}{c} \text{monomials in }\\ \text{ coordinates of } \P_{\Delta}\end{array}\right\} \leftrightarrow \left\{\begin{array}{c} \text{ toric orbifold } \\ \text{ divisors in } \Sigma_{\Delta^{\circ}}\end{array}\right\}.\]

To see why such a correspondence should hold, note that a nonzero element $m \in M \cap \Delta$ can yield two different objects.  On the one hand, since $\P_{\Delta}$ is defined as the closure of the image of the map
\[f: T \rightarrow \P^k\]
\[f(t) = [m_0(t) : \cdots : m_k(t)]\]
in which $\{m_0, \ldots, m_k\} = M \cap \Delta \setminus \{0\}$, the element $m$ gives a coordinate function of $\P_{\Delta}$.  On the other hand, each such $m$ generates a ray in a subdivision of the fan $\Sigma_{\Delta^{\circ}}$ given by cones on faces of $\Delta$.  Thus, $m$ yields a toric divisor in a resolution of $X_{\Sigma_{\Delta^{\circ}}}$, which can be viewed as an orbifold divisor in $X_{\Sigma_{\Delta^{\circ}}}$.  The full monomial-divisor correspondence comes from taking linear combinations of the elements $m \in M \cap \Delta$.

More generally, when $X$ and $X^{\circ}$ are families of Calabi-Yau threefolds, Batyrev's theorem can be phrased in modern terminology as follows:

\begin{theorem}[Batyrev]
\label{Batyrev}
There are isomorphisms
\[H^{1,1}_{CR}(X) \cong H^{2,1}_{CR}(X^{\circ})\]
and
\[H^{2,1}_{CR}(X) \cong H^{1,1}_{CR}(X^{\circ}).\]
\end{theorem}

Note that one must use Chen-Ruan cohomology to account for the presence of orbifold divisors; for a review of Chen-Ruan cohomology, see the Appendix.

To put it another way, one can define a state space for the Calabi-Yau A-model as $H^*_{CR}(X)$, and define a state space for the Calabi-Yau B-model as the same vector space but with a different bi-grading.  Theorem \ref{Batyrev} is then a bi-degree preserving isomorphism between the A-model of $X$ and the B-model of $X^{\circ}$.
 
For smooth threefolds, toric divisors give elements of $H^{1,1}_{CR}(X^{\circ})$, while monomials in the coordinate functions give sections of the anticanonical bundle and these generate $H^{2,1}_{CR}(X)$.  Hence, Theorem \ref{Batyrev} generalizes the monomial-divisor correspondence.

\begin{example}
\label{mirrorquintic}
Let us explore the Batyrev mirror symmetry construction for the quintic threefold $X \subset \P^4$, which is defined by the vanishing of a section of the anticanonical bundle $\O_{\P^4}(5)$.

Generalizing Example \ref{projspace}, the fan for $\P^4$ is generated by the rows of the matrix
\[A:=\left(\begin{array}{cccc} 1 & 0 & 0 & 0 \\ 0 & 1 & 0 & 0\\ 0 & 0 & 1 & 0 \\ 0 & 0 & 0 & 1\\-1 & -1 & -1 & -1  \end{array}\right).\]
If $\Delta$ is the polytope associated to $\P^4$ and its anticanonical divisor, then taking cones on the faces of $\Delta^{\circ}$ should yield the above fan.  From here, it is easy to see that $\Delta^{\circ}$ must be the convex hull of the rows of $A$, and hence the integral points of $\Delta^{\circ}$ are
\[(0,0,0,0), (-1,-1,-1,-1), (1,0,0,0), (0,1,0,0), (0,0,1,0), (0,0,0,1).\]
Thus, $\P_{\Delta^{\circ}}$ is the closure of the image of the morphism
\[f^{\circ}: (\C^*)^4 \rightarrow \P^5\]
\[[z_1:z_2:z_3:z_4] \mapsto [1: z_1^{-1}z_2^{-1}z_3^{-1}z_4^{-1}: z_1 : z_2: z_3: z_4].\]
More explicitly, this shows that $\P_{\Delta^{\circ}} \subset \P^5$ is defined by the equation
\[y_0^5 = y_1y_2y_3y_4y_5.\]

Let $(\Z_5)^3$ act on $\P^4$ diagonally, where we view
\[(\Z_5)^3 = \left\{(\omega_0, \omega_1, \omega_2, \omega_3, \omega_4) \in (\C^*)^5 \; | \; \omega_i^5 =1 \text{ for all } i, \prod_{i=0}^4 \omega_i = 1\right\}/G\]
and $G = \{(\omega, \omega, \omega, \omega, \omega)\}$ is the subgroup of elements that act trivially on $\P^4$.  Then there is a map
\[[\P^4/(\Z_5)^3] \rightarrow \P^5\]
\[[\hat{x}_1: \cdots : \hat{x}_5] \mapsto [\hat{x}_1\hat{x}_2\hat{x}_3\hat{x_4}\hat{x_5}: \hat{x}_1^5: \cdots : \hat{x}_5^5].\]
This map is an isomorphism onto $\P_{\Delta^{\circ}} \subset \P^4$.

The family of anti-canonical hypersurfaces in $\P_{\Delta^{\circ}}$ are given by the restriction of linear functions in the coordinates on $\P^5$.  From the perspective of $[\P^4/(\Z_5)^3]$, these correspond to the $(\Z_5)^3$-invariant quintics
\[a_1 \hat{x}_1^5 + \cdots + a_5\hat{x}_5^5 + a_0 \hat{x}_1 \cdots \hat{x}_5 = 0.\]
Via a change of variables that rescales the coordinates, any such quintic can be expressed as
\[\{x_1^5 + \cdots + x_5^5 + \psi x_1x_2x_3x_4x_5 = 0\} \subset [\P^4/(\Z_5)^3]\]
for a constant $\psi$.  The family of these hypersurfaces (as $\psi$ varies) is the most common expression of the mirror family to the quintic threefold.
\end{example}

\chapter{Hori-Vafa Mirror Symmetry}
\label{HV}

A few definitions are required in order to set the stage for the Hori-Vafa construction.

\section{Basics of symplectic geometry}

References for the material of this section include \cite{Cannas} and \cite{McDuff}.

\begin{definition}
A {\bf symplectic manifold} is a smooth manifold with a closed, nondegenerate $2$-form, referred to as a {\bf symplectic form}.
\end{definition}

Let $M$ be a symplectic manifold, and let $G$ be a Lie group acting on $M$ that preserves the symplectic form $\omega$.  Any element $v$ in the Lie algebra $\mathfrak{g}$ defines a vector field $X_v$ on $M$ giving the infinitesimal action of $v$--- that is,
\[X_v \bigg|_{x} = \left. \frac{d}{dt}\right|_{t=0} \exp(tv) \cdot x\]
for any point $x \in M$, where $\exp: \mathfrak{g} \rightarrow G$ is the exponential map.

\begin{definition}
The action of $G$ on $M$ is {\bf Hamiltonian} if
\begin{enumerate}
\item there exists a {\bf moment map}
\[\mu: M \rightarrow \mathfrak{g}^{\ast},\]
defined by the property that
\[\omega(X_v, \cdot) =  d \bigg( \mu(\cdot)(v) \bigg)\]
for all $v \in \mathfrak{g}$ (here, $\mu(\cdot)(v)$ is a smooth function $M \rightarrow \R$, so its differential is a $1$-form on $M$);
\item the map
\[\mathfrak{g} \rightarrow C^{\infty}(M)\]
\[v \mapsto \mu(\cdot)(v)\]
is a Lie algebra homomorphism, where the Lie bracket on $C^{\infty}(M)$ is defined to be the Poisson bracket.
\end{enumerate}
\end{definition}

Suppose that $(M, \omega)$ is a symplectic manifold equipped with a Hamiltonian action by a group $G$.  Let $s \in \mathfrak{g}^{\ast}$ be a regular value of $\mu$.  Then $\mu^{-1}(s)/G$ has the structure of a symplectic orbifold.  Toric varieties can be constructed as symplectic orbifolds in this way via a method known as {\bf symplectic reduction}, which we describe below.

Start with the symplectic manifold $M= \C^n$, with the standard symplectic structure:
\[\omega = \frac{1}{2} \sum_{i=1}^n dx_i \wedge dy_i = -\frac{1}{2} \text{Im}\left(\sum_{i=1}^n dz_i \wedge d\overline{z_i}\right).\]
Let $G$ be the torus $(S^1)^r$.  An action of $G$ on $M$ is specified by a charge matrix $Q = (Q_{ij})$, where
\[g \cdot x := gQx\]
for a row vector $g \in G$ and a column vector $x \in M$.  The moment map for this action is
\[\mu(z_1, \ldots, z_n) = \left(\frac{1}{2} \sum_{i=1}^n Q_{1i} |z_i|^2, \ldots ,\frac{1}{2} \sum_{i=1}^n Q_{ri} |z_i|^2\right).\]

Any of the symplectic orbifolds obtained as $\mu^{-1}(s)/G$ for a regular value $s$ of this moment map will be toric; namely, they will be of the form
\[X_{\Sigma} = \frac{\C^n \setminus Z(\Sigma)}{(\C^*)^r}\]
for some discriminant locus $Z(\Sigma)$.  The discriminant locus depends on $s$--- it can be thought of as the complement of the image of $\mu^{-1}(s)$ under the action of $(\C^*)^r$--- but only in a rather coarse way.  In particular, $X_{\Sigma}$ is independent of $s$ within ``chambers", connected regions of regular values; only when $s$ crosses a ``wall" at a critical value of $\mu$ will the toric variety change.

\begin{example}
\label{exs}
Consider a toric variety of the form
\[X_{\Sigma} =\frac{\C^{n+1} \setminus Z(\Sigma)}{\C^*},\]
in which the action of $\C^*$ on $\C^{n+1}$ has charge matrix $(1,1, \ldots, 1,-d)$.  Such a variety arises via symplectic reduction on the symplectic manifold $M = \C^{n+1}$, where $G = S^1$ and the moment map is given in coordinates $z_1, \ldots, z_n, p$ on $M$ by
\[\mu = \frac{1}{2}\left(\sum_{i=1}^n |z_i|^2  - d|p|^2\right).\]
The only place where all of the partial derivatives of $\mu$ vanish is $z_1 = \cdots = z_5 = p = 0$, so the only critical value is $s = 0$.  Thus, there are two chambers:

\begin{itemize}
\item If $s > 0$, then the equation $\mu(z_1, \ldots, z_n, p) = s$ implies
\[\sum_{i=1}^n |z_i|^2 = d|p|^2 + \frac{1}{2}s,\]
which can only occur if $\sum_{i=1}^n |z_i|^2 \neq 0$.  The discriminant locus, then, is
\[Z(\Sigma) = \{z_1 = \cdots = z_n = 0\},\]
and we obtain
\[\mu^{-1}(s)/S^1 = \frac{\C^{n+1} \setminus \{z =0\}}{\C^*} = \O_{\P^{n-1}}(-d).\]

\item If $s < 0$, then the equation $\mu(z_1, \ldots, z_n, p) = s$ forces that $p \neq 0$, so
\[Z(\Sigma) = \{p=0\},\]
and
\[\mu^{-1}(s)/S^1 = \frac{\C^{n+1} \setminus \{p=0\}}{\C^*} = [\C^n/\Z_d].\]
\end{itemize}
\end{example}

\section{Gauged Linear Sigma Models}
\label{GLSM}

Using the notion of symplectic reduction, one can define gauged linear sigma models, which are the objects that will be the Hori-Vafa mirrors of varieties.

\begin{definition}
A {\bf gauged linear sigma model}, or {\bf GLSM}, consists of the following data:
\begin{enumerate}
\item an $r \times n$ charge matrix, which can be viewed as the definition of an action of $(\C^*)^r$ on $\C^n$;
\item a moment map $\mu: \C^n \rightarrow \C^r$;
\item an $R$-charge, which is an action of $\C^*$ on $\C^n$ (denoted $\C^*_R$ to distinguish it from the action in (1));
\item a superpotential, which is a map $W: \C^n \rightarrow \C$ satisfying:
\begin{enumerate}
\item $W$ is invariant under the action of $(\C^*)^r$;
\item $W$ is homogeneous of degree $2$ under the action of $\C^*_R$;
\item the critical locus of $W$ is compact inside $\mu^{-1}(s)/(\C^*)^r$ for any regular value $s$ of $\mu$.
\end{enumerate}
\end{enumerate}
\end{definition}

The first two parts of the definition amount to the definition of a toric variety by the technique of symplectic reduction.  The superpotential is constructed so it gives a well-defined map with compact critical locus out of this toric variety.

\begin{example}
\label{glsmex}
The two toric varieties considered in Example \ref{exs} are part of the same GLSM, in which the charge matrix and moment map are as specified above.  The superpotential is
\[W = p(z_1^d + \cdots + z_n^d),\]
which is invariant under the action of $\C^*$ on $\C^{n+1}$ and hence defines a map out of either of the toric varieties constructed above.  Its critical locus in $\mu^{-1}(s)/\C^*$ when $s > 0$ is the hypersurface $\{z_1^d + \cdots + z_n^d = 0\}$ inside the zero-section of $\O_{\P^{n-1}}(-d)$, while its critical locus when $s < 0$ is $\{0\} \subset [\C^n/\Z_d]$.  In either case, this locus is indeed compact.

The $R$-charge of this GLSM can have weights $(0,\ldots, 0,2)$, but this choice is not unique.  For example, another possibility is that the $R$-charge could be $(1,\ldots, 1,-3)$, and these different choices yield possibly different physical theories.
\end{example}

\begin{definition}
A {\bf Landau-Ginzburg model} is a GLSM with a choice of chamber but {\it without} a choice of $R$-charge.  In other words, it consists of a toric variety $X_{\Sigma}$ and a map
\[W: X_{\Sigma} \rightarrow \C\]
whose critical locus is compact.
\end{definition}

\begin{example}
\label{twoLG}
Associated to the GLSM considered in Example \ref{glsmex}, there are two Landau-Ginzburg models, one for each of the chambers.  When $s> 0$, one has the Landau-Ginzburg model
\[W: \O_{\P^{n-1}}(-d) \rightarrow \C\]
\[W = p\sum_{i=1}^5 z_i^d,\]
while for $s<0$, the Landau-Ginzburg model is
\[W: [\C^n/\Z_d] \rightarrow \C\]
\[W = \sum_{i=1}^n z_i^d,\]
since the $p$-coordinate is rescaled to $1$ in this presentation.
\end{example}

The various Landau-Ginzburg models associated to a single GLSM are sometimes called {\bf phases}.

Although a Landau-Ginzburg model need not have an $R$-charge, one can often be chosen, as the above examples illustrate.  A general principle of mirror symmetry is that Calabi-Yau models should have quasihomogeneous mirrors; in other words, the Landau-Ginzburg mirror of a Calabi-Yau hypersurface or complete intersection should have a chosen $R$-charge.

\section{The Jacobian ring}

Hori-Vafa mirror symmetry is an isomorphism of the Calabi-Yau A-model state space associated to a variety (which, recall, is simply its cohomology) with the Landau-Ginzburg B-model state space of a mirror Landau-Ginzburg model.  Specifically, the B-model state space associated to a Landau-Ginzburg model $(Y,W)$, in which $Y$ is a toric variety of dimension $n$, is the {\bf Jacobian ring} (sometimes called the {\bf Milnor ring}):
\[\Jac(W):= \frac{\C[x_1, \ldots, x_n]}{(\d_{x_1}W, \ldots, \d_{x_n}W)}.\]

The situation is especially simple when $W$ is a Morse function.  Indeed, the Jacobian ring can be ``localized" in the sense that
\[\Jac(W) \cong \prod_{x_0} \Jac(W)|_{\text{neighborhood of } x_0},\]
where $x_0$ ranges over critical points of $W$.  If $W$ is Morse, then it takes the form $W = \sum_i y_i^2$ in appropriate coordinates in a neighborhood of each critical point.  It follows that, in the Morse case, the restriction of $\Jac(W)$ to a neighborhood of each critical point is one-dimensional, so the dimension of $\Jac(W)$ equals the number of critical points.  Hori-Vafa mirror symmetry, in this situation, is the claim that if $X$ is a variety with Hori-Vafa mirror $(Y,W)$, then the the number of critical points of $W$ is equal to the dimension of $H^*_{CR}(X)$.

\begin{example}
The $A_n$-singularity is defined as the polynomial $W= x^{n+1}$, viewed as a function $\C \rightarrow \C$.  This has a single critical point at $x=0$, and its Jacobian ring is
\[\Jac(W) = \text{Span}_{\C}\{1, x, x^2, \ldots, x^{n-1}\}.\]
Thus, the dimension of the Jacobian ring equals the number of critical points if and only if $n=1$, which is indeed the only case in which the $A_n$-singularity is Morse.
\end{example}

\begin{example}
\label{hvmirror}
Let $W: (\C^*)^n \rightarrow \C$ be defined by
\[W= x_1 + \cdots x_n + \frac{e^t}{x_1\cdots x_n},\]
where $t$ is an unspecified parameter.  Then
\[\d_{x_i}W = 1 - e^t x_1^{-1} x_2^{-1} \cdots x_n^{-1} x_i^{-1},\]
which vanishes only when
\begin{equation}
\label{critpt}
x_i = e^t x_1^{-1} \cdots x_n^{-1}.
\end{equation}
Since the right-hand side is independent of $i$, there can only be a critical point when $x_1 = \cdots =x_n = \lambda$ for some $\lambda \in \C$.  Equation (\ref{critpt}) implies that $\lambda^{n+1}  = e^t$, so the $n+1$ solutions to this equation will yield $n+1$ critical points.  It is straightforward to check that all of these critical points are nondegenerate, so $W$ is Morse and the dimension of $\Jac(W)$ equals $n+1$.
\end{example}

\section{Hori-Vafa mirrors of compact toric varieties}

We are finally equipped to describe the cohomological statement of Hori-Vafa mirror symmetry.  The material of this section, and of the remainder of the chapter, is based on \cite{HV}.

\begin{definition}
Let
\[X_{\Sigma} = \frac{\C^n \setminus Z(\Sigma)}{(\C^*)^r}\]
be a compact toric variety with charge matrix $Q = (Q_{ij})$.  Then the {\bf Hori-Vafa mirror} is the Landau-Ginzburg model on the toric variety
\[\{x \in \C^n \; | \; x_1^{Q_{1j}} \cdots x_n^{Q_{nj}} = e^{t_j} \text{ for all } 1 \leq j \leq r\}\]
with superpotential given by the restriction of
\[W = x_1 + \cdots + x_n.\]
\end{definition}

One should check that the above subset of $\C^n$ is indeed toric, and that the critical locus of $W$ is compact.  Having done this, the above constitutes a Landau-Ginzburg model, which should be mirror to $X_{\Sigma}$ in the sense that
\[H^*_{CR}(X_{\Sigma}) \cong \Jac(W)\]
whenever $X_{\Sigma}$ is semi-Fano.  Let us verify this in some easy examples.

\begin{example}
The charge matrix of $\P^n$ is $(1,1, \ldots, 1)$, so the underlying toric variety of the Hori-Vafa mirror is the subset of $\C^{n+1}$ defined by the equation
\[x_0 = \frac{e^t}{x_1 \cdots x_n},\]
which is isomorphic to $(\C^*)^n$.  The superpotential is the polynomial considered in Example \ref{hvmirror}.  We saw in that example that $\Jac(W) \cong \C^{n+1}$, so it does match $H^*(\P^n)$ as a vector space.
\end{example}

\begin{example}
A similar computation yields the Hori-Vafa mirror of weighted projective space.  Assume for simplicity that $c_0 = 1$, so $\P(c_0, \ldots, c_n)$ has charge matrix $(1, c_1, \ldots, c_n)$.  Then the superpotential of the Hori-Vafa mirror is
\[W = x_1 + \cdots + x_n + \frac{e^t}{x_1^{c_1} \cdots x_n^{c_n}}.\]
To confirm the statement of mirror symmetry, we compute
\[\d_{x_i}W = 1 - c_ie^tx_1^{-c_1} \cdots x_n^{-c_n} x_i^{-1},\]
so the critical points occur when
\[\frac{x_1}{c_1} = \cdots = \frac{x_n}{c_n} = e^tx_1^{-c_1} \cdots x_n^{-c_n}.\]
If this common value is denoted $\lambda$, then we have
\[\lambda \prod_{i=1}^n (c_i \lambda)^{c_i} = e^t,\]
yielding an equation $\lambda^{1 + c_1 + \cdots + c_n} = Ke^t$ for a constant $K$.  It follows that there are $1 + c_1 + \cdots + c_n$ critical points.  All of these are Morse, so
\[\Jac(W) \cong \C^{1 + c_1 + \cdots + c_n},\]
which is indeed isomorphic to $H^*_{CR}(\P(c_0, c_1, \ldots, c_n))$.
\end{example}

The failure of Hori-Vafa mirror symmetry in the non-semi-Fano case can be observed explicitly for Hirzebruch surfaces $F_n$, which are semi-Fano if and only if $n = 1$.

\begin{example}
Let $F_n$ be the Hirzebruch surface, given by
\[F_n = \frac{\C^4 \setminus Z(\Sigma)}{(\C^*)^2}\]
with charge matrix
\[Q = \left(\begin{array}{cccc} 1 & 1 & n & 0\\ 0 & 0 & 1 & 1\end{array}\right).\]
The constraints defining the toric variety of the mirror Landau-Ginzburg model, then, are
\[x_1x_2x_3^n = e^{t_1},\]
\[x_3x_4 = e^{t_2}.\]
Denoting $q_i = e^{t_i}$ for $i=1,2$, the superpotential is
\[W = x_2 + x_3 + \frac{q_1}{x_2x_3^n} + \frac{q_2}{x_3}.\]
An elementary computation shows that critical points occur only when
\[n^2 q_ 1 x_3^n = (x_3^{n+1} - x_3^{n-1}q_2)^2.\]
This is a polynomial of degree $2n+2$, so for generic values of the parameters $q_1$ and $q_2$, it has $2n+2$ solutions.

On the other hand, $H^*(F_n)$ is $4$-dimensional.  Thus, the number of critical points matches the rank of the cohomology only when $n=1$, which is the semi-Fano case.
\end{example}

\begin{remark}
In recent years, some attempts have been made to adapt Hori-Vafa mirror symmetry to the non-semi-Fano case.  This involves adding higher-order terms to $W$ to ensure that the number of its critical points coincides with the rank of the cohomology of the non-semi-Fano variety.  These correction terms have an interpretation in terms of counting holomorphic discs in $X_{\Sigma}$.

Beyond the state space level, the statements of mirror symmetry at the level of rings and quantum D-modules have also been extended to non-semi-Fano toric varieties, by passing to certain $q$-adic completions.  See the work of Iritani \cite{Iritani} and Gonzalez-Woodward \cite{GW} for details.
\end{remark}

\section{Hori-Vafa mirrors of noncompact toric varieties}

A simple example suffices to show that Hori-Vafa mirror symmetry, in the form stated above, fails when $X_{\Sigma}$ is noncompact:

\begin{example}
Let $X_{\Sigma} = \C$.  Then the Landau-Ginzburg mirror is $W = x: \C \rightarrow \C$.  This has no critical points, so the duality between $H^*(X_{\Sigma})$ and the Jacobian ring fails.
\end{example}

More heuristically, the reason why the above procedure requires compactness involves the aspects of mirror symmetry beyond cohomology.  Roughly speaking, genus-zero mirror symmetry is a correspondence between the quantum cohomology of $X_{\Sigma}$, encapsulated by the $J$-function, and the oscillatory integrals
\[\int_{\Delta} e^{-W/z} d \log(x_1) \cdots d \log(x_n)\]
against cycles $\Delta$.  When $W$ satisfies certain growth properties, the cycles $\Delta$ are in bijection with critical points of $W$, which, in turn, are supposed to correspond via mirror symmetry to the cohomology of $X_{\Sigma}$.  In particular, the $J$-function and the oscillatory integrals depend on the same number of parameters.  If, however, $X_{\Sigma}$ is noncompact, then the $J$-function is not well-defined, and the requisite growth properties of $W$ fail to hold.  Thus, the statement of mirror symmetry breaks down on both sides.

On the other hand, noncompact toric varieties {\it do} have equivariant $J$-functions, so one might expect that there is an equivariant version of the Hori-Vafa mirror for which the mirror symmetry statement is still valid.  This is indeed the case.

\begin{definition}
Let $X_{\Sigma} = (\C^n \setminus Z(\Sigma))/(\C^*)^r$ be a toric variety with charge matrix $Q$.  Then its {\bf equivariant Hori-Vafa mirror} is the Landau-Ginzburg model with superpotential
\[W = x_1 + \cdots + x_n - \sum_{i=1}^n \lambda_i \log(x_i)\]
on the subset of $\C^{n+1}$ defined by the constraints
\[\prod_{i=1}^n x_i^{Q_{ib}} = q_b\]
for nonzero parameters $q_b$.  Here, $\lambda_i$ is a constant, viewed as the equivariant parameter for the $i$th $\C^*$ action.
\end{definition}

\begin{example}
Let $X_{\Sigma} = \C$.  Then the superpotential of the equivariant Hori-Vafa mirror is
\[W = x - \lambda \log(x),\]
so
\[\d _x W = 1 - \frac{\lambda}{x}.\]
Unlike the nonequivariant case, this now has a single critical point, so we recover the correspondence between the dimension of $H^*(X_{\Sigma})$ and the number of critical points.
\end{example}

\begin{example}
\label{Od}
Let $X_{\Sigma}$ be the total space of the bundle $\O_{\P^n}(-d)$, which has charge matrix $(1,1, \ldots, 1, -d)$.  Then
\[W = x_0 + x_1 + \cdots + x_{n+1} - \sum_{i=0}^{n+1} \lambda_i \log(x_i)\]
and the constraint defining the toric variety is
\[x_0 \cdots x_n x_{n+1}^d = q.\]
Thus, we obtain
\[W = x_1 + \cdots + x_{n+1} + \frac{x_{n+1}^d}{x_1 \cdots x_n} - \sum_{i=0}^{n+1} \lambda_i \log(x_i).\]
It is straightforward to check that the number of critical points of this superpotential indeed matches the dimension of $H^*(X_{\Sigma})$.
\end{example}

\begin{example}
It is not always necessary to modify $W$ by all of the terms $\lambda_i \log(x_i)$ in order to achieve mirror symmetry.  For example, consider the toric variety $\O_{\P^1}(-1) \oplus \O_{\P^1}(-1)$, which has charge matrix $Q = (1,1,-1,-1)$.  Its nonequivariant Hori-Vafa mirror is
\[W = x_2 + x_3 + x_4 + q\frac{x_3x_4}{x_2}\]
defined over $(\C^*)^3$, while the equivariant mirror would subtract the terms $\sum_{i=1}^4 \lambda_i \log(x_i)$ from the above.

Even the partial modification
\[\widetilde{W} = x_2 + x_3 + x_4 + q\frac{x_3x_4}{x_2} - \lambda_3 \log(x_3) - \lambda_4 \log(x_4)\]
of $W$ upholds mirror symmetry, though.  This makes sense from the perspective of the $J$-function; only the two $\C^*$ actions in the noncompact fiber directions of $\O_{\P^1}(-1) \oplus \O_{\P^1}(-1)$ are necessary in order to make the equivariant Gromov-Witten theory well-defined.
\end{example}

Despite the failure of nonequivariant Hori-Vafa mirror symmetry for noncompact toric varieties, the noncompact mirrors are still worth remembering, as they play a role in the hypersurface version of mirror symmetry considered below.

\section{The orbifold Jacobian ring}

The Hori-Vafa construction can be adapted to give the mirror of a smooth semi-Fano hypersurface in a toric variety.  As we will see, however, the resulting Landau-Ginzburg model will have a nontrivial symmetry group, and the definition of its B-model state space must be modified accordingly.

Consider a Landau-Ginzburg model of the form $(\C^N, W)$, where $W: \C^N \rightarrow \C$ is a superpotential.  Let $G \subset (\C^*)^N$ be a (finite) group of diagonal matrices preserving $W$; this is referred to as a symmetry group of the Landau-Ginzburg model.

We will define an ``orbifolded" version of the Jacobian ring $\text{Jac}(W)$ that takes the data of $G$ into account.  This is modelled on the definition of the Chen-Ruan cohomology of a global quotient, which is described in the Appendix.  As in the case of Chen-Ruan cohomology, a certain cohomology group is attached to each fixed point of the $G$-action, and the state space is formed by taking the $G$-invariant part of the direct sum of the contributions from all of the fixed points.

\begin{definition}
Given a polynomial $W: \C^N \rightarrow \C$ and a group of symmetries $G$ of the associated Landau-Ginzburg model, the {\bf B-model Landau-Ginzburg cohomology} of $(\C^N,W,G)$ is the vector space
\[\text{Jac}(W,G):= \left(\bigoplus_{g \in G} \text{Jac}(W_g)\right)^G,\]
where
\[W_g=W|_{\text{Fix}(g)}\]
and the $G$-invariant part is taken with respect to the action of $G$ on $\bigoplus_{g \in G}\text{Jac}(W_g)$ that sends
\[\text{Jac}(W_g)\rightarrow \text{Jac}(W_{h^{-1}gh})\]
\[\phi \mapsto \det(h) \cdot h^*\phi.\]
for each $h \in G$.
\end{definition}

This state space will be explained further in Section \ref{LGB} of the next chapter.

\section{Hori-Vafa mirrors of hypersurfaces in toric varieties}

Let us begin by studying the case of Fermat hypersurfaces in projective space.

\begin{example}
Let $X_d \subset \P^{N-1}$ be a smooth degree-$d$ hypersurface defined by the vanishing of the polynomial $A=x_1^d + \cdots +x_N^d$.  Explicitly, the semi-Fano condition corresponds to the requirement that $d \leq N$.

To form the Hori-Vafa mirror, one first constructs a GLSM on a noncompact toric variety for which $X_d$ is the critical locus of the superpotential in a particular phase.  Namely, let $\C^*$ act on $\C^{N+1}$ with charge matrix
\[(1, \ldots, 1, -d),\]
so that the first $N$ factors give precisely the charge matrix for $\P^{N-1}$.  Let the superpotential be
\[W = p \cdot A(x_1, \ldots, x_N): \C^{N+1} \rightarrow \C\]
in coordinates $(x_1, \ldots, x_N, p)$ on $\C^{N+1}$, and let the moment map be
\[\mu = \frac{1}{2} \left(\sum_{i=1}^N |z_i|^2 - d|p|^2\right).\]
(This generalizes the GLSM appearing in Examples \ref{exs} and \ref{glsmex}.)  As we have seen previously, in the $s>0$ phase of this GLSM, the critical locus of $W$ is precisely the hypersurface $X_d$.

Next, construct the (non-equivariant) Hori-Vafa mirror of this noncompact toric variety containing $X_d$.  This is sometimes called the {\bf pre-Hori-Vafa mirror} of $X_d$.  In this case, it is the Landau-Ginzburg model with superpotential
\[\widetilde{W} = x_1 + \cdots + x_N + x_{N+1}\]
on the subset of $\C^{N+1}$ satisfying the constraint
\[x_1 \cdots x_N x_{N+1}^{-d} = e^t.\]
In other words, setting
\[x_i = u_i^d \text{ for } 1 \leq i \leq N\]
and
\[x_{N+1}= u_{N+1},\]
the pre-Hori-Vafa mirror becomes
\[\widetilde{W} = u_1^d + \cdots + u_N^d + e^{-t/d}u_1 \cdots u_N\]
on the toric variety $(\C^*)^N$.

This Landau-Ginzburg model has a nontrivial symmetry group.\footnote{An alternate explanation for the appearance of this automorphism group can be given in terms of the data of the Landau-Ginzburg B-model beyond cohomology.  Specifically, as we have mentioned previously, the full Landau-Ginzburg B-model in genus zero can be viewed as encoding certain oscillatory integrals
\[\int_{\Delta} e^{-W/z} \omega,\]
where $\omega$ is a ``primitive form".  In the toric case, we have mentioned that $\omega = d \log(x_1) \wedge \ldots \wedge d \log(x_n)$.  On the other hand, the primitive form for a hypersurface is $\omega = dx_1 \wedge \ldots \wedge dx_n$.  From this perspective, the appearance of the automorphism group $\text{SL}(W_0)$ is explained by the fact that only automorphisms with determinant $1$ will preserve this form.}  Namely, an automorphism of $(\C^*)^N$ of the form
\[u_i \mapsto \omega_d^{p_i} u_i,\]
for which $\omega_d$ is a $d$th root of unity and
\[\omega_d^{p_1 + \cdots + p_N} = 1,\]
will preserve the superpotential $\widetilde{W}$.  The group of such symmetries is denoted $\text{SL}(W_0)$, since if $W_0 = u_1^d + \cdots + u_N^d$, the automorphisms in question are precisely
\[\left\{\left(\begin{array}{ccc} \lambda_1 && \\ & \ddots & \\ && \lambda_N \end{array}\right) \; \bigg| \; W_0(\lambda_1x_1, \ldots, \lambda_N x_N) = W_0(x_1, \ldots, x_N) \right\} \cap \text{SL}_N(\C).\]

Finally, to form the Hori-Vafa mirror of $X_d$, one ``compactifies" the pre-Hori-Vafa mirror by adding the point $u_1 = \cdots = u_N = 0$ to the domain $\C^N$ of the Landau-Ginzburg model, and then takes the quotient by the above symmetry group.  This yields
\[\widetilde{W}: [\C^N/\text{SL}(W_0)] \rightarrow \C\]
\[\widetilde{W} = u_1^d + \cdots + u_N^d + e^{-t/d} u_1 \cdots u_N\]
as the mirror.  The claim, then, is:
\[H^*_{CR}(X_d) \cong \text{Jac}(\widetilde{W}, \text{SL}(W_0)).\]
That is, the Calabi-Yau A-model cohomology of $X_d$ is isomorphic to the Landau-Ginzburg B-model cohomology of $(\C^N,\widetilde{W}, \text{Aut}(\widetilde{W}))$.

\end{example}

The same basic procedure computes the Hori-Vafa mirror of a hypersurface in a more general weighted projective space.  It is necessary, however, to restrict to a certain class of polynomials.

\begin{definition}
A quasihomogeneous polynomial is {\bf invertible} if the number of monomials equals the number of variables.  That is, after rescaling the variables to absorb any coefficients, the polynomial can be written in the form
\[A = \sum_{i=1}^N \prod_{j=1}^N x_j^{m_{ij}}.\]
\end{definition}

This condition implies that the exponent matrix $E_A = (m_{ij})$ of $A$ is square.  Assuming that $\{A = 0\}$ also defines a smooth orbifold in weighted projective space, it follows moreover that $E_A$ is invertible, which explains the terminology.

\begin{example}
\label{HVex}
Consider a smooth degree-$d$ hypersurface $X$ in weighted projective space $\P(c_1, \ldots, c_N)$ defined by the vanishing of an invertible polynomial $A$ with exponent matrix $(m_{ij})$.  In this case, the semi-Fano condition is $d \leq \sum_{i=1}^N c_i$.

Again, one begins by constructing a GLSM in which $X$ is the critical locus of the superpotential in some phase.  Namely, the GLSM will have charge matrix $(c_1, \ldots, c_N, -d)$ and the superpotential will be $W = p \cdot A(x_1, \ldots, x_N)$.  It follows that the pre-Hori-Vafa mirror of $X$ is the nonequivariant Hori-Vafa mirror of the toric variety $\O_{\P(c_1, \ldots, c_N)}(-d)$ of this GLSM, which is
\[\widetilde{W} = x_1 + \cdots x_N + x_{N+1}\]
defined over the subset of $\C^{N+1}$ with constraint
\begin{equation}
\label{constraint}
x_1^{c_1} \cdots x_N^{c_N} x_{N+1}^{-d} = e^t.
\end{equation}
The constraint can be ``solved"--- that is, expressed in terms of only $N$ variables--- by a change of coordinates:
\begin{equation}
\label{cov}
x_i = \prod_{j=1}^{N} u_j^{m_{ji}}, \;\; i= 1, \ldots, N,
\end{equation}
\[x_{N+1} = u_{N+1}.\]
Then (\ref{constraint}) becomes
\[\prod_{j=1}^N u_j^{\sum_i m_{ji}c_i} u_{N+1}^{-d} = e^t.\]
The fact that $A$  is quasihomogeneous means that
\[\sum_{j=1}^N m_{ij} c_j = d\]
for each $i$, so the constraint is in fact
\[u_{N+1} = e^{-t/d} u_1 \cdots u_N.\]

It follows that the pre-Hori-Vafa mirror has superpotential
\[\widetilde{W} = \sum_{i=1}^N \prod_{j=1}^N u_j^{m_{ji}} + e^{-t/d} u_1 \cdots u_N: (\C^*)^N \rightarrow \C.\]
To form the mirror itself, we add $u_1 = \cdots = u_N = 0$ and take the quotient by $\text{SL}(W_0)$, where $W_0 = \sum_{i=1}^N \prod_{j=1}^N u_j^{m_{ji}}$.  Thus, the mirror is
\[\widetilde{W}: [\C^N/\text{SL}(W_0)] \rightarrow \C\]
with $\widetilde{W}$ as above.
\end{example}

It is interesting to note that the term $W_0$ in the superpotential of the mirror is the transpose of the defining polynomial $A$ of the hypersurface--- that is, the exponent matrix of $W_0$ is the transpose of the exponent matrix for $G$.  The idea that the transpose polynomial should appear in mirror symmetry was actually suggested before Hori and Vafa's work, by physicists Berglund and H\"ubsch, as we will discuss in the next chapter.  At the time when Hori and Vafa proposed their mirror symmetry construction, however, this connection to previous work was not realized.

More generally, the hypersurfaces in toric varieties for which the Hori-Vafa mirror is defined are as follows:

\begin{definition}
\label{inv}
Let $X \subset X_{\Sigma}:= (\C^N \setminus Z(\Sigma))/(\C^*)^r$ be a hypersurface defined by the vanishing of a polynomial $A(x_1, \ldots, x_M)$, where $M \leq N$ and $A$ is homogeneous of degrees $d_1, \ldots, d_r$ with respect to the $r$ actions of $\C^*$.  Let $Q$ be the charge matrix of $ X_{\Sigma}$.  We say that $X$ is {\bf invertible} if
\begin{enumerate}
\item $N-M+1=r$;
\item the $r \times r$ matrix given by taking the last $N-M$ rows of $Q$ and appending the row $(-d_1, \ldots, -d_r)$ is invertible;
\item the exponent matrix of $A$ is square.
\end{enumerate}
\end{definition}

When $X \subset X_{\Sigma}$ is a compact invertible semi-Fano hypersurface, its Hori-Vafa mirror can be constructed by exactly the same procedure as in Example \ref{HVex}.  To summarize, one first forms the Hori-Vafa mirror of the toric variety
\[\frac{(\C^N \setminus Z(\Sigma))\times \C}{(\C^*)^r},\]
in which the charge matrix is given by appending the row $(-d_1, \ldots, -d_r)$ to the charge matrix $Q$ of $X_{\Sigma}$.  This has superpotential
\[\widetilde{W} = x_1 + \cdots + x_{N+1}\]
and is defined on the subset of $\C^{N+1}$ satisfying the constraints
\[\prod_{i=1}^N x_i^{Q_{ij}} x_{N+1}^{-d_j} = e^{t_j}\]
for $j =1, \ldots, r$.  Then, one uses invertibility to express this subset of $
\C^{N+1}$ as $(\C^*)^M$ in variables $u_1, \ldots, u_M$.  Namely, if $(m_{ij})_{i,j =1, \ldots, M}$ is the exponent matrix of $A$ and $(D_{ij})_{i,j = 1, \ldots, r}$ is the invertible matrix described by condition (2) of Definition \ref{inv}, then
\[x_i = \prod_{j=1}^M u_j^{m_{ji}}\]
for $i=1, \ldots, M$ and
\[x_{M+b} = e^{\sum_{j=1}^r t_j D^{jb}} \prod_{i,j,k} u_k^{-m_{ki}Q_{ij}D^{jb}}\]
for $b =1, \ldots, r$, where $(D^{ij})$ is the inverse of $(D_{ij})$.

The Hori-Vafa mirror is then
\[\widetilde{W}: [\C^M/\text{Aut}(\widetilde{W})] \rightarrow \C.\]
Note that the symmetry group $\text{Aut}(\widetilde{W})$ is in general smaller when $N-M$ is larger, since there are more constraints on the variables $u_1, \ldots, u_M$.  This situation occurs when the original polynomial $W = p\cdot A$ has a larger symmetry group.  We will explore the duality between the symmetries of $W$ and those of $\widetilde{W}$ much further in the next chapter.  For now, let us simply look at an example.\footnote{We are grateful to Mark Shoemaker for pointing out this example, and for correcting our definition of invertibility in light of it.}

\begin{example}
Consider the ambient toric variety
\[X_{\Sigma} = \frac{(\C^3 \setminus \{0\}) \times (\C \setminus \{0\})}{(\C^*)^2}\]
with charge matrix
\[Q = \left( \begin{array}{cc} 1 & -2\\ 1 & -1\\ 1 & 0\\ 0 & 3\end{array}\right).\]
It is straightforward to check that $X_{\Sigma} = [\P^2/\Z_3]$.  Let $X$ be the hypersurface defined by the polynomial
\[A = x_1^3 + x_2^3 + x_3^3.\]

To form the Hori-Vafa mirror, first consider the toric variety
\[\frac{(\C^3 \setminus \{0\}) \times (\C \setminus \{0\}) \times \C}{(\C^*)^2}\]
with charge matrix
\[\left( \begin{array}{cc} 1 & -2\\ 1 & -1\\ 1 & 0\\ 0 & 3\\ -3 & 0 \end{array}\right).\]
Its Hori-Vafa mirror (the pre-Hori-Vafa mirror of $X$) is
\[\widetilde{W} = x_1 + x_2  + x_3 + s + p\]
on the subset of $\C^5$ defined by the constraints
\[x_1x_2x_3 p^{-3} = e^{t_1}, \; \;\; x_1^{-2}x_2^{-1}s^3 = e^{t_2}.\]
Setting $x_i = u_i^3$ for $i=1, \ldots, 3$, we find that
\[\widetilde{W} = u_1^3 + u_2^3 + u_3^3 + e^{-t_1/3}u_1u_2u_3 + e^{t_2/3}u_1^2u_2.\]

The symmetry group of $\widetilde{W}$ is not all of $\text{SL}(A^T)$, but instead is isomorphic to $\Z_3$.  Thus, the Hori-Vafa mirror of $X$ is
\[\widetilde{W}: [\C^3/\Z_3] \rightarrow \C\]
for $\widetilde{W}$ as above.
\end{example}

\section[Hori-Vafa mirrors of complete intersections]{Hori-Vafa mirrors of complete intersections in toric varieties}

A similar procedure yields the Hori-Vafa mirror of a complete intersection in a toric variety, but the notion of invertibility must once again be adapted.

\begin{definition}
\label{invci}
Let $X \subset X_{\Sigma} := (\C^N \setminus Z(\Sigma))/(\C^*)^r$ be a complete intersection defined by the vanishing of polynomials
\[A_1(x_1, \ldots, x_M), \ldots, A_k(x_1, \ldots, x_M),\]
where $M \leq N$ and $A_b$ is homogeneous of degree $d_{ba}$ with respect to the $a$th $\C^*$ action.  Let $Q$ be the charge matrix of $X_{\Sigma}$.  We say that $X$ is {\bf invertible} if
\begin{enumerate}
\item $N-M+k =r$;
\item the $r \times r$ matrix given by taking the last $N-M$ rows of $Q$ and appending the matrix $(-d_{ba})$ is invertible;
\item there exists a collection of $N$ monomials
\[\prod_{j=1}^N x_j^{m_{ij}}, \; \; i =1, \ldots, N\]
such that
\[A_b = \sum_{i=1}^N n_{ib} \prod_{j=1}^N x_j^{m_{ij}}\]
with each $n_{ib} \in \{0,1\}$ and $\sum_{b=1}^k n_{ib} = 1$.  (In other words, there are no repeated monomials among the $A_b$, and the total number of monomials appearing is equal to the number of variables.)
\end{enumerate}
\end{definition}

\begin{example}
A complete intersection in $\P^{N_1 -1} \times \cdots \times \P^{N_k -1}$ defined by invertible polynomials $\{A_ i = 0\} \subset \P^{N_i -1}$ is invertible.
\end{example}

\begin{example}
A complete intersection in
\[\P^{N-1} \times \P^{M-1}\]
defined by polynomials
\[A_1 = \sum_{i=1}^N s_i^{d_i}\]
and
\[A_2 = \sum_{j=1}^M s_j t_j^{d_2}\]
with $N \geq M$ is invertible.
\end{example}

As in the case of hypersurfaces, we will proceed by first associating to a compact invertible semi-Fano complete intersection $X \subset X_{\Sigma}$ a GLSM for which $X$ is the critical locus of the superpotential in a particular phase.  From here, a similar set of constraints determined by the defining polynomials of $X$ will yield the pre-Hori-Vafa mirror as a Landau-Ginzburg model.  A partial compactification and quotient by symmetries will convert this into the final Hori-Vafa mirror of $X$.

For simplicity, let us desribe the procedure only in the case where $N=M$, so that $k=r$.  Then $X$ is given by
\[\{A_1 = \cdots = A_k = 0\} \subset \frac{\C^N \setminus Z(\Sigma)}{(\C^*)^k}\]
in which the ambient toric variety $X_{\Sigma}$ has charge matrix $Q = (Q_{ia})$.  Suppose that $A_b$ is homogeneous of degree $d_{ba}$ with respect to the $a$th $\C^*$ action.

The associated GLSM is
\[W: \frac{(\C^N \setminus \{0\}) \times \C^k}{(\C^*)^k} \rightarrow \C,\]
where $(\C^*)^k$ acts with charge matrix
\[\left(\begin{array}{ccc} Q_{11} & \cdots & Q_{1N}\\ \vdots &  & \vdots \\ Q_{k1} & \cdots & Q_{kN}\\-d_{11} & \cdots & -d_{1k}\\ \vdots & & \vdots \\ -d_{k1} & \cdots & -d_{kk}\end{array}\right)\]
and
\[W = p_1A_1(x_1, \ldots, x_N) + \cdots + p_kA_k(x_1, \ldots, x_N)\]
in coordinates $x_1, \ldots, x_N, p_1, \ldots, p_k$ on $\C^N \times \C^k$.

The pre-Hori-Vafa mirror is the Landau-Ginzburg model whose domain is the subset of $\C^{N+1}$ defined by the constraints
\[\prod_{i=1}^N x_i^{Q_{ia}} \prod_{b=1}^k x_N^{-d_{ba}} = e^{t_a}\]
for $a=1, \ldots, a_k$ and whose superpotential is the restriction of
\[\widetilde{W} = x_1 + \cdots + x_N +x_{N+1} + \cdots + x_{N+k}\]
to this subset.  If $m_{ij}$ and $n_{ib}$ are defined as in Definition \ref{invci}, then one can check that the change of variables
\[x_i = \prod_{j=1}^N u_j^{m_{ji}}, \; \; i=1, \ldots, N,\]
\[x_{N+b} = e^{-\sum_{i=1}^k d^{bc} t_c}\prod_{j=1}^N u_j^{n_{jb}}, \;\; b=1, \ldots, k\]
solves the constraint.  Here, $d^{bc}$ is the inverse of the matrix $d_{bc}$ and the expression $-d^{bc}t_c$ uses the Einstein summation convention, and hence should be understood as a sum over $c$.  After this coordinate change, the superpotential of the pre-Hori-Vafa mirror is
\[\widetilde{W} = \sum_{i=1}^N \prod_{j=1}^N u_j^{m_{ji}} + \sum_{b=1}^k e^{-d^{bc}t_c} \prod_{j=1}^N u_j^{n_{jb}}.\]

After adding $u_1 = \cdots = u_N = 0$ and taking the quotient of the domain by symmetries preserving $\widetilde{W}$, one obtains the Hori-Vafa mirror.

\section{An alternative description}

For Calabi-Yau hypersurfaces in toric varieties, we now have two notions of the mirror: the Batyrev mirror and the Hori-Vafa mirror.  Although the two constructions appear unrelated, a different description of the Hori-Vafa mirror can be used to show that it coincides in the Calabi-Yau case with Batyrev's definition.

Recall that for a degree-$d$ semi-Fano hypersurface $X_d = \{A_d = 0\} \subset \P^{N-1}$, the first step in defining the Hori-Vafa mirror is to associate to $X_d$ a GLSM with charge matrix $(1, \ldots, 1, -d)$ and superpotential $W = p \cdot A_d$.  From here, the pre-Hori-Vafa mirror is defined as the Landau-Ginzburg model on the subset of $\C^{N+1}$ satisfying the constraint
\[x_1 \cdots x_N p^{-d} = e^t\]
with superpotential $\widetilde{W} = x_1 + \cdots + x_N + p$.

Previously, we expressed this Landau-Ginzburg model in terms of coordinates $u_1, \ldots, u_N$.  Suppose, however, that we instead used the change of variables
\[p = \tilde{p}\]
\[x_i = \tilde{u_i} \tilde{p} \text{ for } 1 \leq i \leq d,\]
\[x_i = \tilde{u_i} \text{ for } d+ 1 \leq i \leq N.\]
Note that the semi-Fano condition $d \leq N$ is necessary for this to be well-defined.

In these coordinates, the pre-Hori-Vafa mirror is
\[\widetilde{W} = \tilde{p}(\tilde{u_1} + \cdots + \tilde{u}_{d+1}) + \sum_{i=d+1}^N \tilde{u_i}\]
on the subset of $\C^{N+1}$ satisfying the constraint
\[\prod_{i=1}^N \tilde{u_i} = e^t.\]

This GLSM, however, is ``equivalent" to the GLSM with superpotential
\[\widetilde{W} = -\sum_{i=d+1}^N \tilde{u_i}\]
subject to the constraints
\begin{equation}
\label{newconstraint}
\sum_{i=1}^d u_i = -1, \; \; \prod_{i=1}^N u_i = e^t.
\end{equation}
Intuitively, it makes sense that these two theories would be equivalent, since the first of these new constraints kills the first term of the old superpotential.  More explicitly, ``equivalent" means that the two theories yield the same oscillatory integrals.  In particular, since mirror symmetry is defined as a correspondence between the $J$-function and the generating function of such integrals, a theory equivalent to the mirror of $X_d$ will still be mirror.

In the Calabi-Yau case, this new presentation of the mirror agrees with the Batyrev construction.

\begin{example}
The Calabi-Yau condition on $X_d \subset \P^{N-1}$ is $d=N$.  In this case, the new version of the pre-Hori-Vafa mirror described above has no superpotential; it is simply the open manifold defined by (\ref{newconstraint}).  The second equation of (\ref{newconstraint}) implies that the coordinates $\tilde{u_i}$ can be expressed in terms of new variables $z_i$ as
\[\tilde{u_i}= e^{t/N} \frac{z_i^N}{z_1 \cdots z_N},\]
after which the second equation of (\ref{newconstraint}) becomes
\[z_1^N + \cdots + z_N^N + e^{t/N} z_1 \cdots z_N = 0.\]
After the compactification adding $z_1 = \cdots = z_N = 0$ and the quotient by symmetries of $\widetilde{W}$, this is precisely the Batyrev mirror, as we computed in the case $d=N=5$ in Example \ref{mirrorquintic}.
\end{example}

\chapter{Berglund-H\"ubsch-Krawitz Mirror Symmetry}
\label{BHK}

As we have seen, Hori-Vafa mirror symmetry for hypersurfaces reveals an interesting duality between certain polynomials and their transposes: a semi-Fano hypersurface $X = \{W = 0\}$ in weighted projective space is Hori-Vafa mirror to the Landau-Ginzburg orbifold
\begin{equation}
\label{Wtilde}
\widetilde{W} = W^T + e^tx_1 \cdots x_N
\end{equation}
modulo the group $\text{SL}(W^T)$.\footnote{We have changed notation from the previous chapter, denoting the defining polynomial of the hypersurface by $W$ rather than $A$, to be consistent with the literature on LG-to-LG mirror symmetry.}  Given that both $W$ and $W^T$ give rise to Landau-Ginzburg models, this suggests an LG-to-LG version of mirror symmetry.  Our discussion of the resulting statement, known as Berglund-H\"ubsch-Krawitz Mirror Symmetry, is based on \cite{BH} and \cite{Krawitz}.

\section{Phases of the GLSM}

In order to cast Hori-Vafa's Fano-to-LG mirror symmetry in this new framework, we will need to replace the geometric Fano model of the hypersurface $X$ by a Landau-Ginzburg model.  We have already seen how to do this; the trick involves passing from $X$ to an associated gauged linear sigma model.

Suppose $X \subset \P(c_1, \ldots, c_N)$ is a hypersurface in weighted projective space defined by the vanishing of an invertible quasihomogeneous polynomial $W$ of degree $d$.  Then, as observed in the construction of the Hori-Vafa mirror, there is a GLSM for which $X$ is the critical locus in one phase.  Namely, let
\[\overline{W} := p \cdot W : \frac{\C^N \times \C}{\C^*} \rightarrow \C,\]
where the charge matrix for the action of $\C^*$ is $(c_1, \ldots, c_N, -d)$.  The moment map for this GLSM is
\[\mu = \frac{1}{2}\left(\sum_{i=1}^N c_i |x_i|^2 - d|p|^2\right).\]
The only critical value is $s=0$, so the GLSM has two phases.  The important fact for Hori-Vafa mirror symmetry is that $X$ is the critical locus of $\overline{W}$ in the phase $s>0$, which is the Landau-Ginzburg model
\[\overline{W}: \O_{\P(c_1, \ldots, c_N)}(-d) \rightarrow \C,\]
as computed in Example \ref{twoLG} of the previous chapter.  This is sometimes referred to as the ``geometric phase", since it arises out of a hypersurface.  The $s<0$ phase, sometimes called the ``Landau-Ginzburg phase", is
\[W: [\C^N/\Z_d] \rightarrow \C,\]
where the generator of $\Z_d$ acts on $\C^N$ via the matrix
\begin{equation}
\label{J}
J := \left(\begin{array}{ccc} e^{2\pi i \frac{c_1}{d}} && \\ & \ddots & \\ & & e^{2\pi i \frac{c_N}{d}} \end{array}\right) \in U(1).
\end{equation}

A Landau-Ginzburg model admits two types of cohomology.  The one that we have seen thus far, the orbifold Jacobian ring of the superpotential, is the {\bf Landau-Ginzburg B-model cohomology}.  On the other hand, there is a {\bf Landau-Ginzburg A-model cohomology}, which we will define later in this chapter.

When applied to the geometric phase of the above GLSM, the Landau-Ginzburg A-model cohomology is simply the cohomology of the hypersurface $X$.  Thus, the Hori-Vafa mirror symmetry statement
\[H^*_{CR}(X) \cong \Jac(\widetilde{W}, \text{SL}(W^T))\]
can be phrased as an exchange of the Landau-Ginzburg A-model cohomology of the $s> 0$ phase of the GLSM with the Landau-Ginzburg B-model cohomology of the model $(\C^N, \widetilde{W}, \text{SL}(W^T))$ defined in the previous chapter.

To be more precise, the Hori-Vafa mirror $\widetilde{W}$ of $X$ depended on a parameter $t$, and in fact, $X$ is only mirror to this Landau-Ginzburg model when $t$ is near $\infty$--- this is the {\bf large complex structure limit} of the family of Landau-Ginzburg models.  On the other hand, if $t$ is near $-\infty$, then the B-model Landau-Ginzburg cohomology of $(\C^N, \widetilde{W}, \text{SL}(W^T))$ will instead match the A-model of the $s<0$ phase of the GLSM.\footnote{Roughly, $t$ is the real part of the parameter $s$ in the GLSM.}

In the special case where $X$ is Calabi-Yau, the Landau-Ginzburg B-models of $\widetilde{W}$ for different values of the parameter $t$ can be related to one another.  We thus obtain a diagram:
\begin{equation}
\label{diag1}
\xymatrix{
\text{LG A-model of } \overline{W}|_{s>0}\ar@{<->}[d]_{\text{Hori-Vafa}} & \text{LG A-model of } \overline{W}|_{s<0}\ar@{<->}[d]^{\text{Hori-Vafa}}\\
\text{LG B-model of } \widetilde{W}|_{t \rightarrow \infty}\ar@{<-->}[r] & \text{LG B-model of } \widetilde{W}|_{t \rightarrow -\infty},
}
\end{equation}
where each Landau-Ginzburg model should be understood as orbifolded with respect to its symmetry group.  The dotted horizontal arrow is a special case of a conjecture known as the Landau-Ginzburg/Calabi-Yau (LG/CY) correspondence, which relates the various phases of a GLSM to one another.  The conjecture has been proven in many instances; see \cite{CR} \cite{CIR} \cite{Acosta} \cite{ClRo}, among many other references, for more information.

Given the discussion of the previous paragraphs and the fact that $\widetilde{W} \rightarrow W^T$ as $t \rightarrow \infty$, the diagram (\ref{diag1}) can be re-expressed as:
\begin{equation}
\label{diag2}
\xymatrix{
	\text{CY A-model of }X\ar@{<->}[d] & \text{LG A-model of }(W,\Z_d)\ar@{<->}[d]\\
	\text{LG B-model of } (W^T, \text{SL}(W^T))\ar@{<-->}[r]\ar@{<-->}[ur] &  \text{LG B-model of } (\widetilde{W}|_{t \rightarrow -\infty}, \text{Aut}(\widetilde{W})).
}
\end{equation}
The arrow from the bottom-left to the top-right, which is implied by the LG/CY correspondence, is Berglund-H\"ubsch-Krawitz mirror symmetry.

More generally, Berglund-H\"ubsch-Krawitz mirror symmetry can be extended by replacing the group $\Z_d$ with any group $G$ of diagonal symmetries of  $\C^N$ that preserves $W$ and contains the element $J$ defined in (\ref{J}).  The construction yields a ``mirror group" $G^T$ associated to any such $G$.  The generalization of the diagonal arrow in (\ref{diag2}), then, is
\[\text{LG A-model of }(W,G) \cong \text{LG B-model of } (W^T, G^T),\]
the latter of which is Hori-Vafa mirror to the orbifold $[X/(G/\langle J \rangle)]$.

All three forms of mirror symmetry discussed in these notes can now be schematically related.  Let $X_W:=\{W=0\}$ be a Calabi-Yau hypersurface.  For a group $G$ of symmetries as above, let $\widetilde{G}$ denote $G/\langle J \rangle$.  Then we have a diagram:
\[
\xymatrix{
	\text{CY A-model of } [X_W/\widetilde{G}]\ar@{<->}[d]_{\text{Batyrev}}\ar@{<-->}[r]\ar@{<->}[dr]^{\text{HV}} & \text{LG A-model of } (W,G)\ar@{<->}[d]^{\text{BHK}}\ar@{<->}[dl]\\
	\text{CY B-model of } [X_{W^T}/\widetilde{G^T}]\ar@{<-->}[r] & \text{LG B-model of } (W^T,G^T),
}
\]
in which the horizontal arrows are again the LG/CY correspondence.

It should be noted that, although we have used the Calabi-Yau assumption to motivate the appearance of Berglund-H\"ubsch-Krawitz mirror symmetry, the state space isomorphism between $(W,G)$ and $(W^T, G^T)$ holds even when the Calabi-Yau assumption fails.  The remainder of this chapter will be devoted to making the specific assumptions and results of Berglund-H\"ubsch-Krawitz mirror symmetry precise.

\section{Classification of nondegenerate singularities}

The polynomial $W$ will be required to be invertible, quasihomogeneous, and to satisfy a certain nondegeneracy condition.

\begin{definition}
A quasihomogeneous polynomial $W$ is {\bf nondegenerate} if
\begin{enumerate}
\item zero is the only critical point;
\item there are no monomials of the form $x_i x_j$ with $i \leq j$.
\end{enumerate}
\end{definition}

A few remarks about this definition are in order.  First, the condition that zero is the only critical point is equivalent to the existence of only isolated critical points, since quasihomogeneity implies that $W$ and its derivatives vanish at $\mathbf{x}$ if and only if the same is true for any scalar multiple $\lambda \mathbf{x}$. 

Second, suppose that we define:

\begin{definition}
The {\bf exponent matrix} of an invertible quasihomogeneous polynomial
\[W = \sum_{i=1}^N a_i \prod_{j=1}^N x_j^{m_{ij}}\]
is the matrix $E_W:=(m_{ij})$.
\end{definition}

\begin{definition}
The {\bf charges} of a quasihomogeneous polynomial with weights $c_1, \ldots, c_N$ and degree $d$ are the rational numbers
\[q_j = \frac{c_j}{d}.\]
\end{definition}

Then the definition of nondegeneracy has the following elementary consequences:
\begin{enumerate}
\item The charges $q_j$ are unique.
\item $q_j \leq 1/2$ for all $j$.
\item The exponent matrix of $W$ is nonsingular.
\end{enumerate}

Polynomials satisfying the above conditions admit a very nice classificaiton.

\begin{theorem}[Kreuzer-Skarke \cite{KS}] Suppose that $W$ is a nondegenerate quasihomogeneous invertible polynomial.  Then $W$ can be written as a sum of polynomials
\[W = \sum_{s=1}^k W_s\]
in disjoint sets of variables, where each $W_s$ is of one of the following types:
\begin{enumerate}
\item Fermat: $x^a$ for some $a \geq 2$;
\item Loop: $x_1^{a_1}x_2 + x_2^{a_2}x_3 + \cdots + x_N^{a_N}x_1$;
\item Chain: $x_1^{a_1}x_2 + x_2^{a_2}x_3 + \cdots x_N^{a_N}$.
\end{enumerate}
\end{theorem}

This is extremely useful for computations in the Landau-Ginzburg model, because the model decomposes as a ``product", in a precise sense, whenever the superpotential breaks up as a disjoint sum.  As a result, the study of the LG model reduces to the study of the model associated to Fermat, loop, and chain polynomials separately.  Moreover, the operation of transpose preserves the types, so mirror symmetry can also be studied in the three cases individually.

Let us mention a few specific examples of nondegenerate invertible quasihomogeneous polynomials, following \cite{Arnold}.  An important invariant of such a polynomial is its {\bf central charge}
\begin{equation}
\label{centralcharge}
c_W := \sum_{i=1}^N (1 - 2q_i).
\end{equation}
When $\sum q_i = 1$--- or, equivalently, $c_W = N-2$--- the polynomial is said to be {\bf Calabi-Yau}, since it defines a Calabi-Yau hypersurface in weighted projective space.

The only polynomials with $c_W < 1$ are the {\bf ADE}-singularities
\begin{itemize}
\item $A_n = \frac{1}{n+1} x^{n+1}$;
\item $D_n = x^{n-1} + xy^2, n \geq 4$;
\item $E_6 = x^3 + y^3$;
\item $E_7 = x^3 + xy^3$;
\item $E_8 = x^3 + y^5$.
\end{itemize}
These all have $\sum q_i < 1$, so none of them is Calabi-Yau.  Note, furthermore, that $A$-type and $E$-type polynomials are self-mirror under the operation of transpose, whereas $D$-type polynomials are not self-mirror.

Among the Calabi-Yau examples, those with $c_W = 1$ are the {\bf elliptic singularities}
\begin{itemize}
\item $P_8 = x^3 + y^3 + z^3$,
\item $X_9 = x^2 + y^4 + z^4$,
\item $J_{10} = x^2 + y^3 + z^6$.
\end{itemize}
The Calabi-Yau examples with $c_W= 2$ have four variables.  These are the {\bf K3 singularities}, of which there are $95$ types.  When $c_W = 3$, we have the Calabi-Yau threefolds, for which there are thousands of examples.  As one can see, classification becomes unwieldy beyond this point.

\section{The maximal diagonal symmetry group}

Let $W$ be a nondegenerate quasihomogeneous invertible polynomial.  The {\bf maximal diagonal symmetry group} of $W$ is the group of diagonal matrices preserving $W$; that is,
\[G_{max} = \left\{g= \left( \begin{array}{ccc} e^{2\pi i g_1} && \\ & \ddots & \\ && e^{2\pi i g_N}\end{array}\right) \; \bigg| \; W(g\cdot \mathbf{x}) = W(\mathbf{x}) \right\}.\]
The groups $G$ for which the Landau-Ginzburg A-model and B-model are defined will always be subgroups of $G_{max}$.

There is a convenient description of this group in terms of the exponent matrix of $W$.  By the definition of invertibility, $W$ is of the form
\[W = \sum_{i=1}^N \prod_{j=1}^N x_j^{m_{ij}},\]
for a matrix $E_W = (m_{ij})$.  Let
\[\rho_k = \left(\begin{array}{c} \rho_1^{(k)} \\ \vdots \\ \rho_N^{(k)} \end{array}\right), \; \; k = 1, \ldots, N\]
be the columns of the inverse matrix $E_W^{-1}$.  Convert each into a diagonal matrix
 \[\overline{\rho_k} := \left( \begin{array}{ccc} e^{2\pi i \rho_1^{(k)}} && \\ & \ddots & \\ && e^{2\pi i \rho_N^{(k)}}\end{array}\right).\]

\begin{proposition}
The maximal diagonal symmetry group $G_{max}$ is generated by the matrices $\overline{\rho_1}, \ldots, \overline{\rho_N}$.
\begin{proof}
The first thing to check is that each $\overline{\rho_k}$ lies in $G_{max}$.  Indeed,
\begin{align*}
W(\overline{\rho_k} \mathbf{x}) & = W(e^{2\pi i \rho_1^{(k)}}x_1, \ldots, e^{2\pi i \rho_N^{(k)}} x_N)\\
&=\sum_{i=1}^N\prod_{j=1}^N e^{2\pi i \rho_j^{(k)}m_{ij}} x_j^{m_{ij}}\\
&=\sum_{i=1}^N \left(e^{2\pi i \sum_{j} \rho_j^{(k)}m_{ij}}\right) \prod_{j=1}^N x_j^{m_{ij}}.
\end{align*}
The exponent $\sum_j \rho_j^{(k)}m_{ij}$ is precisely the $i$th entry in the vector $E_W \cdot \rho_k$, which is the $(i,j)$th entry in the matrix $E_W \cdot E_W^{-1}$.  Hence, it equals $\delta_{ij}$, so we obtain
\[W(\overline{\rho_k} \mathbf{x}) = \sum_{i=1}^N \prod_{j=1}^N x_j^{m_{ij}} = W(\mathbf{x}),\]
as required.

Now, to see that $\overline{\rho_1}, \ldots, \overline{\rho_N}$ generate $G_{max}$, notice that for a diagonal matrix
\begin{equation}
\label{g}
g= \left( \begin{array}{ccc} e^{2\pi i g_1} && \\ & \ddots & \\ && e^{2\pi i g_N}\end{array}\right) \in G_{max},
\end{equation}
we have
\[E_W  \left(\begin{array}{c} g_1 \\ \vdots \\ g_N \end{array}\right) = g_1 E_W\rho_1 + \cdots + g_NE_W\rho_N = E_W(g_1 \rho_1 + \cdots + g_N \rho_N),\]
since the column vectors $\rho_k$ are defined so that $E_W \rho_k$ is a vector with a $1$ in the $k$th entry and zeroes elsewhere.  Cancelling $E_W$ from both sides of (\ref{g}) shows that
\[g = g_1 \overline{\rho_1} + \cdots + g_N \overline{\rho_N},\]
so $g$ can be expressed in terms of $\overline{\rho_1}, \ldots, \overline{\rho_N}$.
\end{proof}
\end{proposition}

Note, in particular, that the distinguished element $J \in G_{max}$ defined in ({\ref{J}) is presented in these generators as
\[J = \overline{\rho_1} \cdots \overline{\rho_k},\]
since it is the diagonal matrix corresponding to the column vector $(q_1, \ldots, q_N)^T$ and
\[E_W \cdot  \left(\begin{array}{c} q_1 \\ \vdots \\ q_N \end{array}\right) =  \left(\begin{array}{c} 1 \\ \vdots \\ 1 \end{array}\right).\]

\section{The mirror group}

Given a quasihomogeneous invertible polynomial $W$, denote its maximal diagonal symmetry group by $G_{W, max}$, and denote the maximal diagonal symmetry group of the transpose polynomial by $G_{W^T, max}$; recall, $W^T$ is defined by the relationship
\[E_{W^T} = (E_W)^T.\]
For any subgroup $G \subset G_{W,max}$ containing $J$, we will construct a ``mirror" subgroup $G^T \subset G_{W^T,max}$.  Before doing so, two simple observations are useful.

First, whereas $G_{max, W}$ consists of matrices $\rho = \rho_1^{k_1} \cdots \rho_N^{k_N}$, where $\rho_i$ are the diagonal matrices associated to columns of $W$, the group $G_{max, W^T}$, consists of matrices $h = h_1^{r_1} \cdots h_N^{r_N}$, in which $h_i$ are associated to the {\it rows} of $W$.

Second, the same proof that showed that $\rho_k \in G_{max}$ implies that a matrix
\[g = \left( \begin{array}{ccc} e^{2\pi i g_1} && \\ & \ddots & \\ && e^{2\pi i g_N}\end{array}\right)\]
lies in the maximal diagonal symmetry group $G_{max}$ of $W$ if and only if
\[E_W^{-1} \left(\begin{array}{c} g_1 \\ \vdots \\ g_N \end{array}\right) \in \Z^N.\]
A similar observation holds for the maximal diagonal symmetry group of $W^T$, but replacing right multiplication by a column with left multiplication by a row.

With these facts in mind, we define $G^T$ as follows.

\begin{definition}
Given a subgroup $G \subset G_{W,max}$ containing $J$, the {\bf mirror (or transpose) group} is
\[\left\{h_1^{r_1} \cdots h_N^{r_N} \; \bigg| \; \begin{array}{ccc}(r_1 & \cdots & r_N)\end{array} E_W^{-1}  \left(\begin{array}{c} k_1 \\ \vdots \\ k_N \end{array}\right) \in \Z^N \text{ for any } \rho_1^{k_1} \cdots \rho_N^{k_N} \in G\right\}.\]
\end{definition}

To put it another way, one can define a pairing between $G_{max, W}$ and $G_{max, W^T}$ by
\[\langle \rho, h \rangle = \begin{array}{ccc}(r_1 & \cdots & r_N)\end{array} E_W^{-1}  \left(\begin{array}{c} k_1 \\ \vdots \\ k_N \end{array}\right) \mod \Z,\]
and in terms of this pairing, $G^T$ is the orthogonal complement of $G$.

The following are some useful observations regarding the operation of transpose on both polynomials and groups:
\begin{enumerate}
\item If $G_1 \subset G_2$, then $G_2^T \subset G_1^T$.
\item $(G^T)^T = G$.
\item $\{1\}^T = G_{max}$.
\item $\langle J \rangle^T = \text{SL}({W^T})$.
\end{enumerate}
In particular, these facts imply that the transpose of groups preserves the {\bf Calabi-Yau condition}
\[\langle J \rangle \subset G \subset \text{SL}({W})\]
for subgroups $G$.  This is important for aspects of LG-to-LG mirror symmetry that will not be discussed in these notes; namely, there are ring structures on the A-model and B-model LG state spaces, and mirror symmetry gives a ring isomorphism only in the Calabi-Yau case.

\section{B-model LG cohomology}
\label{LGB}

We have now defined the mirror $(W^T, G^T)$ of a Landau-Ginzburg model $(W,G)$.  In order to make the statement of Berglund-H\"ubsch-Krawtiz mirror symmetry precise, we must carefully define the cohomology groups on the A-side and the B-side that will be exchanged.

The Landau-Ginzburg B-model cohomology has already appeared in the context of Hori-Vafa mirror symmetry, where we viewed it as an orbifold Jacobian ring.  Let us recall and expand upon this definition.

It is convenient to give a somewhat different presentation of the Jacobian ring.  Let
\[\Omega_W = \Omega^N(\C^N)/\big(dW \wedge \Omega^{N-1} (\C^N)\big).\]
Then the map
\[\Jac(W) \rightarrow \Omega_W\]
\[\phi \mapsto \phi \cdot dx_1 \cdots dx_n =: \phi d \underline{x}\]
is an isomorphism of vector spaces.\footnote{Some care should be taken here, because both sides admit both a grading and a $G$-action, and the isomorphism does not preserve these aspects.  We will be careful to specify these in what follows, but the basic principle is that they are inherited from $\Omega_W$, not from $\Jac(W)$.}

There is a pairing on $\Omega_W$ defined by
\begin{align*}
\langle \phi d\underline{x}, \phi' d \underline{x} \rangle &= \text{Res}\left(\frac{\phi \phi' d\underline{x}}{\d_1W \cdots \d_NW}\right)\\
&= \frac{1}{2\pi i} \int_{|\d_i W| = \epsilon} \frac{\phi \phi' d \underline{x}}{\d_1 W \cdots \d_NW}.
\end{align*}

Furthermore, there is a grading
\[\deg(\phi d\underline{x}) = \deg(\phi) + \sum_{i=1}^N q_i\]
on $W$, in which the {\bf charges} $q_i$ are defined by
\[q_i = \frac{c_i}{d}.\]
Under this grading, the pairing can equivalently be computed by using the fact that, like the cohomology of a manifold, $\Omega_W$ has a unique generator in the top degree.  This top degree is the central charge $c_W$ defined in (\ref{centralcharge}), and the generator is the Hessian $\text{Hess}(W)$.  Just as one computes the pairing on the cohomology of a manifold via the volume form, then, one can compute
\[\langle \phi, \phi'\rangle = \lambda,\]
where
\[\phi \cdot \phi' = \frac{\lambda}{\mu} \text{Hess}(W) + \text{ lower-order terms}.\]
Here, the normalizing constant is the Milnor number $\mu$, defined by
\[\mu = \dim(\Jac(W)) = \prod_{i=1}^N \left(\frac{1}{q_i} - 1 \right).\]

The B-model Landau-Ginzburg cohomology can now be built by orbifolding, as was previously described.

\begin{definition}
Given a nondegenerate quasihomogeneous invertible polynomial $W$ and a subgroup $G \subset G_{W, max}$ containing $J$, the {\bf B-model Landau-Ginzburg cohomology} of $(W,G)$ is the vector space
\begin{equation}
\label{sectors}
\Omega_{W,G} = \left(\bigoplus_{g \in G} \Omega_{W_g}\right)^G,
\end{equation}
where
\[W_g = W|_{\text{Fix}(g)}.\]
Here, the $G$-invariant part is taken with respect to the action of $G$ on $\bigoplus_{g \in G} \Omega_{W_g}$ that sends
\[\Omega_{W_g} \rightarrow \Omega_{W_{h^{-1}gh}}\]
via pullback under multiplication by $h$.
\end{definition}

We should note that the restriction of $W$ to $\text{Fix}(g)$ is still nondegenerate, so each $\Omega_{W_g}$ is of the form described previously.

The pairings on the vector spaces $\Omega_{W_g}$ defined above can be combined to give a pairing on $\Omega_{W,G}$, again in much the same way that the Poincar\'e pairing is defined in Chen-Ruan cohomology.  Specifically, one pairs
\[\Omega_{W_g}^G \otimes \Omega_{W_{g^{-1}}}^G \rightarrow \C\]
with the previously-defined residue pairing, using the fact that $\text{Fix}(g) = \text{Fix}(g^{-1})$.

As in Chen-Ruan cohomology, a shift in the grading on $\Omega_{W,G}$ is necessary in order to make the pairing behave like a Poincar\'e pairing.  First, define an unshifted (or ``internal") bigrading on $\Omega_W$ by doubling the single grading:
\[\Omega^{p,q}_{W_g} = \begin{cases} (\Omega_{W_g})^p & \text{ if } p = q\\ 0 & \text{ if } p \neq q.\end{cases}\]
Then, define a {\bf degree shift}:
\[(Q_-^B, Q_+^B) = (\iota_{(g^{-1})}, \iota_{(g)}) - \left(\sum_{i=1}^Nq_i, \sum_{i=1}^Nq_i\right).\]
Here, $\iota_g$ is the age shift from Chen-Ruan cohomology, defined by
\[\iota_{(g)} = \sum_{i=1}^N \frac{m_i}{m},\]
where $m = \text{ord}(g)$ and $g$ acts on the tangent space $T_xX$ to a point $x \in \text{Fix}(g)$ by
\[g = \left(\begin{array}{ccc} e^{2\pi i \frac{m_1}{m}} && \\ & \ddots & \\ && e^{2\pi i \frac{m_n}{m}}\end{array}\right).\]
The grading on $\Omega_{W,G}$ is then defined as
\[\Omega^{p,q}_{W,G} = \left(\bigoplus_{g \in G} \Omega_{W_g}^{p - Q_-^B, q - Q_+^B}\right)^G.\]

The purpose of this shift is to ensure that the following result holds:

\begin{lemma}
\label{degs}
If
\[\langle \phi_1 d\underline{x}_g, \phi_2 d \underline{x}_{g^{-1}} \rangle \neq 0,\]
then
\[\deg_W^{\pm}(\phi_1 d \underline{x}_g) + \deg_W^{\pm}(\phi_2 d \underline{x}_{g^{-1}}) = c_W.\]
\begin{proof}
We will require the fact that, if
\[g = \left(\begin{array}{ccc} e^{2\pi i \Theta^1_g} && \\ & \ddots & \\ && e^{2\pi i \Theta^N_g}\end{array}\right),\]
then
\[\iota_{(g)} = \sum_{\Theta^i_g \neq 0} \Theta^i_g,\]
whereas
\[\iota_{(g^{-1})} = \sum_{\Theta^i_g \neq 0} (1- \Theta^i_g).\]
Using this, we have
\begin{align*}
&\deg^{\pm}_W(\phi_1 d\underline{x}_g) + \deg^{\pm}_W(\phi_2 d \underline{x}_{g^{-1}})\\
=& \deg(\phi_1 d \underline{x}_g) + Q^B_{\pm}(g) + \deg(\phi_2 d\underline{x}_{g^{-1}}) + Q^B_{\pm}(g^{-1})\\
=&\deg(\phi_1) + \sum_{\Theta^i_g = 0} q_i + Q^B_{\pm}(g) + \deg(\phi_2) + \sum_{\Theta^i_g = 0} q_i + Q^B_{\pm} (g^{-1})\\
=&\deg(\phi_1) + \deg(\phi_2) + \sum_{\Theta^i_g = 0}(1-2q_i) + 2 \sum_{\Theta^i_g = 0} q_i + \sum_{\Theta^i_g \neq 0} 1 - 2\sum_{i=1}^Nq_i\\
=&c_{W_g} + \sum_{\Theta^i_g \neq 0} (1 - 2q_i)\\
=&\sum_{i=1}^N(1- 2q_i)\\
=& c_W.
\end{align*}
\end{proof}
\end{lemma}

The upshot of Lemma \ref{degs} is that $\Omega_{W,G}$ behaves like the cohomology of a complex manifold of dimension $c_W$--- despite the fact that $c_W$ may be fractional.  This is referred to as a ``manifold of dimension $c_W$" in the physics literature to make sense of the notion of fractional dimension.

The components of the decomposition (\ref{sectors}) are referred to as {\bf sectors}.  There are a number of special cases that yield particularly important sectors.

First, when $g=1$, the component of $\Omega_{W,G}$ is called the {\bf nontwisted sector}.  It is isomorphic to $\Omega_{W}$ but is graded via the shift
\[(Q^B_-, Q^B_+) = -\left(\sum_{i=1}^Nq_i, \sum_{i=1}^Nq_i\right).\]
Thus,
\[\deg^{\pm}_W(\phi d\underline{x}) = \deg(\phi).\]
The element $1 d \underline{x}$ of bidegree $(0,0)$ and the element $\text{Hess}(W) d \underline{x}$ of bidegree $(c_W, c_W)$ both lie in the nontwisted sector.

The $J$-sector, for which $g = J$, is $1$-dimensional:
\[\Omega^G_{W_J} \cong \C,\]
since $\text{Fix}(J) = \{0\}$.  It degree shift is
\[(Q^B_-, Q^B_+) = \left(\sum_{i=1}^N (1-q_i), \sum_{i=1}^Nq_i \right) - \left(\sum_{i=1}^Nq_i, \sum_{i=1}^Nq_i\right)= (c_W, 0).\]
The $J^{-1}$-sector is also one-dimensional, but its degree shift is
\[(Q^B_-, Q^B_+) = (0,c_W),\]
by an analogous computation.

Thus, the components of various bidegrees in the Landau-Ginzburg B-model cohomology can be compiled into the following rough outline of a Hodge diamond:
\[
\xymatrix{
& H^{c_W, c_W} & \\
\C \cong \Omega_{W_J} \cong H^{c_W, 0} & & H^{0, c_W} \cong \Omega_{W_{J^{-1}}} \cong \C\\
 & H^{0,0}, &
}
\]
where the middle column is the nontwisted sector.

\section{A-model LG cohomology}

We have actually mentioned the A-model Landau-Ginzburg cohomology before briefly, as well.  Analogously to the way in which the B-model cohomology is built out of Jacobian rings, the building blocks of the A-model cohomology are the vector spaces
\begin{equation}
\label{middlecohom}
H^N(\C^N, W^{+ \infty}; \C),
\end{equation}
where
\[W^{+ \infty} = (\text{Re}W)^{-1}(\rho, \infty)\]
for $\rho \gg 0$.\footnote{Somtimes $W^{+ \infty}$ is replaced with a ``Milnor fiber" $W^{-1}(z_0)$ of $W$, where $z_0$ is a sufficiently large real number.  This does not affect the vector space; however, in other parts of the theory it is necessary to use the Hodge structure on cohomology, and in order for a natural Hodge structure to be defined, it is necessary to use $W^{+ \infty}$ and not a single fiber.}

Just as in the B-model case, we will define the A-model by taking a direct sum of middle cohomology groups as in (\ref{middlecohom}) for restrictions of $W$ to the fixed-point sets of $g \in G$.  There will be a degree shift, defined so that a certain pairing behaves like a Poincar\'e pairing.  Towards this end, we will need a pairing on (\ref{middlecohom}).

The first step in that direction is to define a pairing
\begin{equation}
\label{Apairing}
H^N(\C^N, W^{-\infty}; \Q) \otimes H^N(\C^N, W^{+\infty};\Q) \rightarrow \Q,
\end{equation}
where
\[W^{-\infty} = (\text{Re}W)^{-1}(-\infty, -\rho)\]
for $\rho \gg 0$.  To do so, we use the fact that $H^N(\C^N, W^{+ \infty}, \Q)$ is dual to the homology group $H_N(\C^N, W^{+ \infty};\Q)$.  The latter has a basis consisting of the preimages in $\C^N$ of a collection of nonintersecting paths in $\C$ that begin at the critical values of $W$ and move in the direction of $\text{Re}(z) = \infty$, eventually becoming horizontal.  These subsets of $\C^N$ are called {\bf Lefschetz thimbles}; see Figure \ref{thimbles}.  There is an analogous basis for $H_N(\C^N, W^{-\infty};\C)$, consisting of ``opposite" thimbles.

\begin{figure}
\begin{tikzpicture}
\draw (0,0) -- (5,0);
\draw (5,0) -- (6,2);
\draw (6,2) -- (1,2);
\draw (1,2) -- (0,0);
\draw[dashed] (3.5,0) -- (4.5,2);
\draw ( 1.85,0.5) circle[radius=2pt];
\fill ( 1.85,0.5) circle[radius=2pt];
\draw (1.85, 0.5) -- (5.25, 0.5);

\draw (1.5, 0.95) circle[radius=2pt];
\fill (1.5, 0.95) circle[radius=2pt];
\draw (1.5, 0.95) -- (5.475, 0.95);

\draw (2.5, 1.4) circle[radius=2pt];
\fill (2.5, 1.4) circle[radius=2pt];
\draw (2.5, 1.4) -- (5.7, 1.4);

\draw (1, 0.25) circle[radius=2pt];
\fill (1, 0.25) circle[radius=2pt];
\draw plot [smooth, tension=1] coordinates {(1,0.25) (0.85, 1) (2.25,1.7)};
\draw (2.25, 1.7) -- (5.85, 1.7);

\draw (5.75, 1) node[right] {\C};
\end{tikzpicture}
\caption{The preimages in $\C^N$ of these paths are closed at one end and open at the other, giving them the appearance of infinite ``thimbles".}
\label{thimbles}
\end{figure}
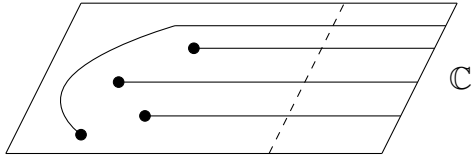

The two types of thimbles can intersect one another, so there is an intersection pairing
\[H_N(\C^N, W^{-\infty};\Q) \otimes H_N(\C^N, W^{+\infty}; \Q) \rightarrow \Q.\]
In fact, it is a perfect pairing, and (\ref{Apairing}) is obtained by dualization.  To be more explicit, if $\delta_i$ is a basis for $H_N(\C^N, W^{-\infty};\Q)$ and $\delta_i^{\vee}$ is a basis for $H_N(\C^N, W^{+\infty}; \Q)$, then
\[\langle \alpha, \beta \rangle = \sum_{i,j} \langle \alpha, \delta_i \rangle \langle \delta_i, \delta_j^{\vee} \rangle \langle \beta, \delta_j^{\vee} \rangle,\]
where the three pairings on the right-hand side are given by intersection in $H_N(\C^N, W^{+ \infty};\Q)$, $H_N(\C^N;\Q)$, and $H_N(\C^N, W^{- \infty};\Q)$, respectively.

From here, (\ref{Apairing}) can be converted into a pairing on (\ref{middlecohom}) by using the morphism
\[I: \C^N \rightarrow \C^N\]
\[(x_1, \ldots, x_N) \mapsto (\xi^{c_1}x_1, \ldots, \xi^{c_N}x_N)\]
for a chosen $\xi$ satisfying $\xi^d = -1$.  This morphism has
\[W(I(x_1, \ldots, x_N)) = -W(x_1, \ldots, x_N),\]
and hence it interchanges $W^{+ \infty}$ with $W^{-\infty}$.  The pairing on (\ref{middlecohom}) is defined as
\[\langle \alpha, \beta \rangle :=  \langle \alpha, I^*\beta \rangle.\]

Although the resulting pairing appears to depend on a choice of $\xi$, we will always look at its restriction to
\begin{equation}
\label{smaller}
H^N(\C^N, W^{+ \infty} ;Q)^{G}
\end{equation}
for $G$ a subgroup of $G_{max}$ containing $J$.  One can check that the pairing on $H^N(\C^N, W^{+\infty}; \Q)^{\langle J \rangle}$ is well-defined and independent of $\xi$, so the same is true for the restriction to the smaller spaces (\ref{smaller}).

The orbifolding construction is very similar to what we have seen previously.

\begin{definition}
If $W$ is a nondegenerate quasihomogeneous invertible polynomial and $G \subset G_{W,max}$ is a subgroup containing $J$, then the {\bf A-model Landau-Ginzburg cohomology} of $(W,G)$ is the vector space
\begin{equation}
\label{Adecomp}
\mathcal{H}_{W,G}:= \left(\bigoplus_{g \in G} H^{N_g}(\C^N_g, W^{+\infty}_g; \Q)\right)^G,
\end{equation}
where $\C^N_g$ is the fixed locus of $g$, $W^{+\infty}_g = \text{Re}(W_g^{-1}(\rho, \infty))$ for $\rho \gg 0$, $W_g = W|_{\C^N_g}$, and $N_g$ is the complex dimension of $\C^N_g$.  The action of $G$ takes the sector indexed by $g$ to the sector indexed by $h^{-1}gh$ via pullback under multiplication by $h$.
\end{definition}

This has a pairing, given by mapping
\[H^{N_g}(\C^N_g, W^{+ \infty}_g; \Q)^G \otimes H^{N_{g^{-1}}}(\C^N_{g^{-1}}, W^{+ \infty}_{g^{-1}}; \Q)^G \rightarrow \Q\]
via the pairing defined previously; note that $\C^N_g = \C^N_{g^{-1}}$.  There is also a degree shift, defined so that this pairing behaves like a Poincar\'e pairing.  It is:
\[\mathcal{H}^{p,q}_{W,G} = \left(\bigoplus_{g \in G} \Omega_{W_g}^{p - Q_-^A, q - Q_+^A}\right)^G,\]
where the shift is 
\[(Q_-^A, Q_+^A) := (\iota_g, \iota_g) - \left(\sum_{i=1}^N q_i, \sum_{i=1}^N q_i\right)\]
and the internal bi-grading is given by the Hodge structure on the vector space $H^{N_g}(\C^N_g, W^{+ \infty}_g; \Q)$.

In fact, the A-model and B-model Landau-Ginzburg cohomology of $(W,G)$ are isomorphic to one another, but with different gradings.  To prove this, one shows that the map
\[\Omega_W \rightarrow H^N(\C^N, W^{+ \infty};\C) = \Hom(H_N(\C^N, W^{+ \infty}; \Z), \C)\]
\[\phi d \underline{x} \mapsto \left( \Delta \mapsto \int_{\Delta} e^{- W} \phi d \underline{x} \right)\]
is an isomorphism; it also respects the pairing on either side, up to a constant.\footnote{This is a more computationally practical way to compute the A-model pairing than via Lefschetz thimbles.}  It does not respect the grading, though; instead, $\Omega^p_W$ maps to the bidegree $(N-p,p)$ part of $H^N(\C^N, W^{+ \infty};\C)$.

Even after the degree shift, the two sides are bigraded differently.  An element of $(\bigoplus_g \Omega_{W_g})^G$ with internal bidegree $(p,p)$ and degree shift
\[(Q^-_B, Q^+_B) = (\iota_{g^{-1}}, \iota_g) - \left(\sum q_i, \sum q_i\right)\]
corresponds to an element in $\left(\bigoplus_g H^{N_g}(\C^N_g, W^{+\infty}_g;\Q)\right)^G$ with internal degree $(N_g - p,p)$ and degree shift
\[(Q^-_A, Q^+_A) = (\iota_g, \iota_g) - \left( \sum q_i, \sum q_i\right).\]
One can compute, then, that
\[\deg^+_A = \deg^+_B\]
but
\[\deg^- A = c_W - \deg^-_B.\]
In particular, this says that the $A$-model Landau-Ginzburg cohomology of $(W,G)$ is obtained from the $B$-model cohomology of $(W,G)$ by flipping the Hodge-diamond.

Let us use this relationship to the B-model cohomology to study a few of the special sectors in the decomposition (\ref{Adecomp}).  The nontwisted sector, indexed by the element $g=1$, has degree shift
\[(Q_-^A, Q_+^A) = \left(-\sum_{i=1}^Nq_i, -\sum_{i=1}^N q_i\right)\]
and internal grading
\[\deg(\phi d\underline{x}) = \deg(\phi) + \sum_{i=1}^N q_i.\]
In particular, the A-model image of the element $1 d \underline{x}$ has bidegree
\begin{align*}
&\big(\deg^-_A(d\underline{x}), \deg^-A(d\underline{x})\big)\\
= &(N- \sum q_i, \sum q_i) + (-\sum q_i, -\sum q_i)\\
= &(c_W, 0),
\end{align*}
and the A-model image of the element $\text{Hess}(W) d \underline{x}$ has bidegree
\begin{align*}
&\big(\deg^-_A(\text{Hess}(W)d\underline{x}), \deg^-_A(\text{Hess}(W)d \underline{x})\big)\\
=&(N- c_W - \sum q_i, c_W + \sum q_i)+ (-\sum q_i, - \sum q_i)\\
=&(0, c_W).
\end{align*}

The sector $g = J$ has $\C^N_J = \{0\}$, so the middle-dimensional relative cohomology is generated by the constant function $1$, which we denote $e_J$.  This has
\[\big(\deg^-_A(e_J), \deg^+_A(e_J)\big) = \left(\sum q_i, \sum q_i\right) + \left(-\sum q_i, - \sum q_i\right) = (0,0).\]
The sector $g = J^{-1}$ also has a single generator $e_{J^{-1}}$, with bidegree $(c_W, c_W)$.  Comparing with the sketch of the B-model Hodge diamond computed in the previous section, we again see that passing between the A-grading and the B-grading corresponds to a rotation.

\section{Berglund-H\"ubsch-Krawitz mirror symmetry}

The result that makes the association $(W,G) \mapsto (W^T, G^T)$ qualify as mirror symmetry is the following:

\begin{theorem}[Krawitz \cite{Krawitz}]
\label{BHKthm}
There is a bidegree-preserving isomorphsim
\[\Omega_{W,G} \cong \mathcal{H}_{W^T, G^T}.\]
\end{theorem}

Alternatively, the results of the previous section show that this can be expressed as an isomorphism between the B-model of $(W,G)$ and the B-model of $(W^T, G^T)$ that rotates the bigrading.

For $g \in G$, let $F_g \subset \{1, \ldots, N\}$ be the indices of the coordinate directions fixed by $g$, so that
\[\text{Fix}(g) = \text{Spec} \left( \C\left[\left\{ x_i \; | \; i \in F_g\right\}\right] \right).\]
An element of $\Omega_{W,G}$ will be written as
\[\sum_{i \in F_g} x_i^{r_i} dx_i \left| \prod_{j=1}^N \rho_j^{s_j + 1}\right\rangle\]
if it is drawn from the sector indexed by
\[g = \prod_{j=1}^N \rho_j^{s_j + 1} \in G.\]
Under this notation, the mirror map in Theorem \ref{BHKthm} is given in terms of the B-model Landau-Ginzburg cohomology by
\[\Omega_{W,G} \rightarrow \Omega_{W^T, G^T}\]
\[\prod_{i \in F_g} x_i^{r_i} d x_i \left| \prod_{j=1}^N \rho_j^{s_j+1}\right\rangle \mapsto \prod_{j=1}^N y_j^{s_j} dy_j \left| \prod_{i \in F_g} \rho_i^{r_i +1}\right\rangle.\]
That is, it exchanges group elements with monomials in the coordinates.

\begin{example}
Let
\[W = x^3y + xy^5,\]
for which $W^T = W$.  This has weights $c_x = 2$ and $c_y = 1$ and degree $7$, so
\[J = \left(\begin{array}{cc} e^{2\pi i\frac{2}{7}} & 0 \\ 0 & e^{2\pi i\frac{1}{7}} \end{array}\right).\]
Let
\[ G = \langle J\rangle,\]
for which 
\[G^T = \text{SL}({W^T}) = \left\langle \left( \begin{array}{cc}  e^{2\pi i\frac{1}{2}} & 0 \\ 0 & e^{2\pi i\frac{1}{2}} \end{array}\right)\right\rangle.\]

By computing the inverse of the exponent matrix, it is easy to verify that the generators of $G_{\text{max}, W}$ are
\[\rho_x = \left(\begin{array}{cc} e^{2\pi i\frac{5}{14}} & 0 \\ 0 & e^{2\pi i\frac{-1}{14}} \end{array}\right)\]
and
\[\rho_y = \left(\begin{array}{cc} e^{2\pi i\frac{-1}{14}} & 0 \\ 0 & e^{2\pi i\frac{3}{14}} \end{array}\right).\]
There are seven elements in the subgroup $\langle J \rangle$.  Using the fact that $J = \rho_x \rho_y$ and the identity $\rho_x^3 \rho_y=1$, these can be written as
\[\langle J \rangle = \{\rho_x^0\rho_y^0, \rho_x^1\rho_y^1, \rho_x^2\rho_y^2, \rho_x^3\rho_y^3, \rho_x^1\rho_y^3, \rho_x^2\rho_y^4, \rho_x^3\rho_y^5\}.\]
All but the first of these is narrow, so the corresponding sector of $\Omega_{W,\langle J \rangle}$ will just be
\[\Q \cdot e_{\rho_x^a \rho_y^b},\]
in which $e_{\rho_x^a \rho_y^b}$ denotes the volume form on the sector indexed by $\rho_x^a \rho_y^b$.  There is no need to restrict to the $\langle J \rangle$-fixed part on these sectors, since the action is trivial.  For the remaining sector, indexed by $\rho_x^0\rho_y^0$, the action of $J$ sends
\[x^i y^j dx dy \mapsto \text{det}(J) \cdot (e^{2\pi i\frac{2}{7}}x)^i (e^{2\pi i \frac{1}{7}}y)^j dx dy = e^{2\pi i\frac{3 + 2i + j}{7}} x^a y^b dx dy.\]
As a result, the $\langle J \rangle$-invariant part is spanned by
\[\{x^2 dxdy, xy^2 dxdy, y^4 dxdy\}.\]
(Since we are working in the Jacobian ring
\[\Jac(W) = \frac{\C[x,y]}{(3x^2y + y^5, x^3 + 5xy^4)},\]
powers of $y$ greater than $4$ can always be expressed in terms of smaller exponenents.)

As for the mirror side, there are two elements in $\langle J\rangle^T$, which can be written as
\[\langle J \rangle^T = \{\rho_x^0\rho_y^0, \rho_x^2 \rho_y^3\}.\]
The second of these gives the sector
\[\Q \cdot e_{\rho_x^2 \rho_y^3}\]
of $\Omega_{W^T, \langle J \rangle^T}$.  The first gives the nontwisted sector.  Since the nontrivial element of $\langle  J \rangle^T$ acts by
\[x^i y^j dx dy \mapsto \text{det}(\rho_x^2\rho_y^3) \cdot (e^{2\pi i\frac{1}{2}}x)^i (e^{2\pi i \frac{1}{2}}y)^j dx dy = e^{2\pi i\frac{i + j}{2}} x^a y^b dx dy,\]
the $\langle J \rangle^T$-invariant part of the nontwisted sector is spanned by
\[\{1, xy, x^2y^2, y^2, xy^3, x^2y^4, x^2, y^4\},\]
again keeping the powers of $y$ below $5$.

These generators are matched up by the mirror map in the following way:

\begin{center}
\begin{tabular}{|c||c|c|c|c|}
\hline
$(W,J)$ & $e_{\rho_x^1\rho_y^1}$ & $e_{\rho_x^2\rho_y^2}$ &$ e_{\rho_x^3\rho_y^3}$ &$ e_{\rho_x^1\rho_y^3}$\\
\hline
$(W^T,\langle J \rangle^T)$ & $1e_{\rho_x^0\rho_y^0}$ & $xy e_{\rho_x^0\rho_y^0}$ & $x^2y^2 e_{\rho_x^0\rho_y^0}$ & $y^2 e_{\rho_x^0\rho_y^0}$\\
\hline
\end{tabular}
\end{center}

\vspace{0.25cm}

\begin{center}
\begin{tabular}{|c||c|c|c|c|c|}
\hline
$(W,J)$ & $e_{\rho_x^2\rho_y^4}$ &$ e_{\rho_x^3\rho_y^5}$ & $x^2 e_{\rho_x^3 \rho_y^1}$ &  $xy^2 e_{\rho_x^0}$ & $y^4 e_{\rho_x^1\rho_y^5}$\\
\hline
$(W^T, \langle J \rangle^T)$ & $xy^3 e_{\rho_x^0\rho_y^0}$ & $x^2y^4 e_{\rho_x^0\rho_y^0}$ & $x^2 e_{\rho_x^3 \rho_y^1}$ & $e_{\rho_x^2\rho_y^3}$ & $y^4 e_{\rho_x^1\rho_y^5}$\\
\hline
\end{tabular}
\end{center}

\end{example}

Notice that the mirror map sends ``narrow sectors"--- those indexed by group elements $g$ for which $\text{Fix}(g) = \{0\}$--- to elements of the nontwisted sector.  In the Calabi-Yau case, there are exactly as many narrow sectors for $(W,G)$ as elements of the nontwisted sector for $(W^T, G^T)$, and the matching is perfect.

Via the Landau-Ginzburg/Calabi-Yau correspondence, this sets up a strong parallel between the Berglund-H\"ubsch-Krawitz and Batyrev-Borisov mirror symmetry constructions.  Indeed, in all known cases of the LG/CY correspondence, the narrow sectors of the Landau-Ginzburg cohomology correspond to the ambient part of the cohomology of $X_W$, consisting of classes pulled back from the weighted projective space.  Such ambient casses are generated by toric divisors.  Thus, the exchange between narrow group elements and nontwisted monomials in Berglund-H\"ubsch-Krawitz mirror symmetry is matched, via the correspondence, with the exchange between toric divisors and monomials in Batyrev mirror symmetry.

\chapter*{Appendix: Chen-Ruan Cohomology}

Many more details on this topic can be found in \cite{ALR}.

An orbifold, speaking geometrically, is an object that is locally the quotient of a manifold by the action of a finite group.  One can make this definition precise in the category-theoretic language of groupoids.  From this perspective, an orbifold is a groupoid for which the objects and arrows form manifolds and the structure morphisms (source, target, composition, identity, and inversion) are all smooth.

The case on which we will focus is when $X$ is a complex manifold and $G$ is a finite group acting on $X$.  Then there is an orbifold 
\[\mathcal{X} = [X/G].\]
This should be thought of as a version of the quotient that records any isotropy of the original action.  The fact that the orbifold ``remembers" the data of the action is clear from the groupoid point of view: the orbifold groupoid $[X/G]$ has objects $X$ and arrows $x \rightarrow g \cdot x$ for each $x \in X$ and $g \in G$.

The de Rham or singular cohomology of an orbifold can be defined, but it is insufficient for capturing all of the data of $\mathcal{X}$.  One way to understand the problem is through Gromov-Witten theory.  Ordinarily, in Gromov-Witten theory, there would be evaluation morphisms
\[\ev_i: \M_{g,n}(Y, \beta) \rightarrow Y\]
to record the images of the various marked points.  In the orbifold setting, however, a morphism $f: \mathcal{C} \rightarrow \mathcal{X}$ from an {\it orbifold} curve $\mathcal{C}$ to an {\it orbifold} $\mathcal{X}$ has more local data around a marked point $x_i$ than simply its image.  Namely, $\mathcal{C}$ is of the form $[\C/G_i]$ near $x_i$ for a finite group $G_i$, and part of the data of $f$ is a homomorphism $G_i \rightarrow G$.  Thus, the evaluation maps should keep track not only of the images of the marked points but of the homomorphisms on isotropy, and for this reason, they should land not in $\mathcal{X}$ but in a more complicated object.

\begin{definition}
The {\bf inertia orbifold} of $\mathcal{X} = [X/G]$ is
\[\Lambda [X/G] = \left[ \bigsqcup_{g \in G} (X^g \times \{g\})/G\right],\]
where
\[X^g = \{x \in X\; | \; gx = x\},\]
and $G$ acts on this disjoint union by
\[h\cdot (x,g) = (hx, hgh^{-1}).\]
\end{definition}

Equivalently, one can write
\[\Lambda [X/G] = \left[ \bigsqcup_{(g) \in G^*} (X^g \times \{g\})/C(g)\right],\]
in which $G^*$ denotes the set of conjugacy classes of $G$.  The advantage of this presentation is that it breaks $\Lambda[X/G]$ into its connected components.  In particular, there is a component
\[[(X \times \{1\})/G] \subset \Lambda[X/G]\]
that is isomorphic to $\mathcal{X}$ itself, called the {\bf nontwisted sector}.  The other connected components are referred to as {\bf twisted sectors}.

The inertia orbifold is the image of the evalaution maps in Gromov-Witten theory, so, at least from that perspective, its cohomology is the natural object to study.

\begin{definition}
As a vector space, the {\bf Chen-Ruan cohomology} of $\mathcal{X}$ is defined as the cohomology of $\Lambda[X/G]$.
\end{definition}

However, the grading on Chen-Ruan cohomology differs from that of the inertia stack.

\begin{definition}
The {\bf degree (or age) shift} of an element $g \in G$ is defined as
\[\iota(g):= \sum_{i=1}^n \frac{m_{i,g}}{m_g},\]
where $m_g = \text{ord}(g)$ and $g$ acts on the tangent space $T_xX$ to a point $x \in X^g$ by
\begin{equation}
\label{action}
\rho_x(g) := \left(\begin{array}{ccc} e^{2\pi i \frac{m_{1,g}}{m_g}} && \\ & \ddots & \\ && e^{2\pi i \frac{m_{n,g}}{m_g}}\end{array}\right).
\end{equation}
\end{definition}

One can check that $\iota$ gives a locally constant function on $\Lambda[X/G]$.  The grading on Chen-Ruan cohomology is defined as follows:
\[H^{p,q}_{CR}([X/G]; \Q) := \prod_{(g) \in G^*} H^{p- \iota(g), q - \iota(g)} (X^g/C(g); \Q).\]

There is an involution $I: X^g \times \{g\} \rightarrow X^{g^{-1}} \times \{g^{-1}\}$ for any $g \in G$, which is simply the identity on the first component.  Using this, one can define a pairing on Chen-Ruan cohomology by
\[\langle \; , \; \rangle : H^{p,q}(X^g \times \{g\}/C(g); \Q) \otimes H^{p', q'}(X^{g^{-1}} \times \{g^{-1}\}/C(g^{-1}); \Q) \rightarrow \Q.\]

\begin{proposition}
Let $n = \dim_{\C}(X)$, and choose $\alpha \in H^{p,q}(X^g\times \{g\}/C(g))$ and $\beta \in H^{p', q'}(X^{g^{-1}} \times \{g^{-1}\}/C(g^{-1}))$.  If $\langle \alpha, \beta \rangle \neq 0$, then
\[p + \iota(g) + p' + \iota(g^{-1}) = n,\]
\[q + \iota(g) + q' + \iota(g^{-1}) = n.\]
That is, the degrees after shifting must sum to $n$.
\begin{proof}
We claim that
\[\iota(g) + \iota(g^{-1}) = \text{rank}(\rho_x(g) - I),\]
with $\rho_x(g)$ as in (\ref{action}).  Indeed, $m_{i,g^{-1}} \equiv -m_{i,g} \mod m_g$, so if $m_{i,g} \neq 0$, then $m_{i,g} + m_{i,g^{-1}} = m_{i,g}$ and hence
\[\frac{m_{i,g}}{m_g} + \frac{m_{i,g^{-1}}}{m_g} = 1.\]
If $m_{i,g} = 0$, however, then $m_{i,g^{-1}} = 0$, and so
\[\frac{m_{i,g}}{m_g} + \frac{m_{i,g^{-1}}}{m_g} = 0.\]
It follows that $\iota(g) + \iota(g^{-1})$ is the number of entries not equal to $1$ in $\rho_x(g)$, which is precisely the rank of $\rho_x(g) - I$.

The dimension of the fixed locus in $T_xX$ under the action of $g$ is the same as the dimension of $X^g$, so the above implies that
\[\iota(g) + \iota(g^{-1}) = n - \dim_{\C}(X^g).\]
Thus, the Proposition follows from the fact that two elements of ordinary cohomology pair nontrivially only if
\[p + p' = q + q' = \dim_{\C}(X^g).\]
\end{proof}
\end{proposition}

The role of the degree shift is now clear: it makes the pairing on Chen-Ruan cohomoloyg behave with respect to degrees like the Poincar\'e pairing on the ordinary cohomology of a manifold.  There are further parallels between Chen-Ruan and ordinary cohomology.  For example, there is a product structure, with a unit lying in the nontwisted sector whose Poincar\'e dual is the volume form, also drawn from the nontwisted sector.  However, we only need to study the Chen-Ruan cohomology as a graded vector space in these notes.

\begin{example}
Let $\Z_r$ act on $\C$ via multiplication by $r$th roots of unity.  Then
\[\Lambda [\C/\Z_r] = [\C/\Z_r] \sqcup [\bullet/\Z_r] \sqcup \cdots \sqcup [\bullet/\Z_r],\]
in which there are $r-1$ twisted sectors in addition to the nontwisted sector.  Thus, as a vector space, the Chen-Ruan cohomology of $[\C/\Z_r]$ is $\Q^r$, with generators $e_0, e_1, \ldots, e_{r-1}$ given by the constant function $1$ on each of the above components.  The degree-shift on the component generated by $e_i$ is $\frac{i}{r}$, so the $r$ copies of $\Q$ occur in bidegrees $(0,0), (\frac{1}{r}, \frac{1}{r}), \ldots, (\frac{r-1}{r}, \frac{r-1}{r})$.
\end{example}

\bibliographystyle{abbrv}
\bibliography{BModelBib}

\begin{thebibliography}{10}

\bibitem{Acosta}
P.~Acosta.
\newblock Asymptotic expanesion and the {LG}/({F}ano, {G}eneral {T}ype)
  correspondence.
\newblock {\em arXiv preprint arXiv:1411.4162}, 2014.

\bibitem{ALR}
A.~Adem, J.~Leida, and Y.~Ruan.
\newblock {\em Orbifolds and stringy topology}, volume 171.
\newblock Cambridge University Press, 2007.

\bibitem{Arnold}
V.~I. Arnolʹd.
\newblock {\em Singularity theory}, volume~53.
\newblock Cambridge University Press, 1981.

\bibitem{Batyrev}
V.~Batyrev.
\newblock Dual polyhedra and mirror symmetry for {C}alabi-{Y}au hypersurfaces.
\newblock {\em J. Algebraic Geom.}, 3:493--535, 1994.

\bibitem{BB}
V.~Batyrev and L.~Borisov.
\newblock On {C}alabi-{Y}au complete intersections in toric varieties.
\newblock In M.~Andreatta and T.~Peternell, editors, {\em Higher-dimensional
  complex varieties (Trento, 1994)}. Walter de Gruyter, Berlin, 1996.

\bibitem{BH}
P.~Berglund and T.~H\"ubsch.
\newblock A generalized construction of mirror manifolds.
\newblock {\em Nuclear Physics B}, 393, 1994.

\bibitem{Candelas}
P.~Candelas, X.~C. de~la Ossa, P.~Green, and L.~Parkes.
\newblock A pair of {C}alabi-{Y}au manifolds as an exactly soluble
  superconformal field theory.
\newblock {\em Nuclear Phy. B}, 359:21--74, 1991.

\bibitem{Cannas}
A.~Cannas~da Silva.
\newblock {\em Lectures on Symplectic Geometry}.
\newblock Number 1764 in Lecture Notes in Mathematics. Springer-Verlag, 2008.

\bibitem{CaK}
E.~Cattani and A.~Kaplan.
\newblock Degenearting variations of {H}odge structures.
\newblock In {\em Actes du {C}olloque du {T}h\'eorie de {H}odge, {L}uminy
  1987}, volume 179--180, pages 67--96. Soc. Math. France, 1989.

\bibitem{CR}
W.~Chen and Y.~Ruan.
\newblock A new cohomology theory of orbifold.
\newblock {\em Communications in Mathematical Physics}, 248(1):1--31, 2004.

\bibitem{CIR}
A.~Chiodo, H.~Iritani, and Y.~Ruan.
\newblock Landau-{G}inzburg/{C}alabi-{Y}au correspondence, global mirror
  symmetry and {O}rlov equivalence.
\newblock {\em Publications math{\'e}matiques de l'IH{\'E}S}, pages 1--90,
  2012.

\bibitem{CR2}
A.~Chiodo and Y.~Ruan.
\newblock Landau--{G}inzburg/{C}alabi--{Y}au correspondence for quintic
  three-folds via symplectic transformations.
\newblock {\em Inventiones mathematicae}, 182(1):117--165, 2010.

\bibitem{CR1}
A.~Chiodo and Y.~Ruan.
\newblock {LG/CY} correspondence: the state space isomorphism.
\newblock {\em Advances in Mathematics}, 227(6):2157--2188, 2011.

\bibitem{Clader}
E.~Clader.
\newblock Landau-{G}inzburg/{C}alabi-{Y}au correspondence for the complete
  intersections ${X}_{3,3}$ and ${X}_{2,2,2,2}$.
\newblock {\em arXiv preprint arXiv:1301.5530}, 2013.

\bibitem{ClRo}
E.~Clader and D.~Ross.
\newblock Sigma models and phase transitions for complete intersections.
\newblock {\em In preparation}, 2015.

\bibitem{CL}
E.~A. Coddington and N.~Levinson.
\newblock {\em Theory of ordinary differential equations}.
\newblock McGraw-Hill, 1955.

\bibitem{CK}
D.~Cox and S.~Katz.
\newblock {\em Mirror symmetry and algebraic geometry}, volume~68 of {\em
  Mathematical {S}urveys and {M}onographs}.
\newblock Amer. Math. Soc., 1999.

\bibitem{CLS}
D.~A. Cox, J.~B. Little, and H.~K. Schenck.
\newblock {\em Toric varieties}.
\newblock American Mathematical Society, 2011.

\bibitem{Deligne1}
P.~Deligne.
\newblock {\em Equation Diff\'erentielles \'a Points Singularier r\'eguliers},
  volume 163 of {\em Lectures Notes in Mathematics}.
\newblock Springer-Verlag, 1970.

\bibitem{FJR}
H.~Fan, T.~Jarvis, and Y.~Ruan.
\newblock The {W}itten equation, mirror symmetry and quantum singularity
  theory.
\newblock {\em Annals of Mathematics}, 178:1--106, 2013.

\bibitem{Fulton}
W.~Fulton.
\newblock {\em Introduction to Toric Varieties}.
\newblock Princeton University Press, 1993.

\bibitem{Givental1}
A.~Givental.
\newblock Equivariant {G}romov-{W}itten invariants.
\newblock {\em Internat. Math. Res. Notices}, 13:613--663, 1996.

\bibitem{Givental2}
A.~Givental.
\newblock A mirror theorem for toric complete intersections.
\newblock In {\em Topological field theory, primitive forms, and related
  topics}, volume 160 of {\em Progr. Math.}, pages 141--175. Birkh\"auser,
  1998.

\bibitem{GW}
E.~Gonzalez and C.~Woodward.
\newblock Quantum cohomology and toric minimal model programs.
\newblock {\em arXiv preprint arXiv:1207.3253}, 2012.

\bibitem{Griffiths1}
P.~Griffiths.
\newblock On the periods of certain rational integrals: {I}, {II}.
\newblock {\em Annals of Mathematics}, pages 460--495, 498--541, 1969.

\bibitem{Griffiths2}
P.~Griffiths, editor.
\newblock {\em Topics in transcendental algebraic geometry}.
\newblock Number 106 in Ann. of Math. Stud. Princeton University Press, 1984.

\bibitem{Hori}
K.~Hori, S.~Katz, R.~Pandharipande, R.~Thomas, C.~Vafa, R.~Vakil, and
  E.~Zaslow.
\newblock {\em Mirror Symmetry}.
\newblock American Mathematical Society/ Clay Mathematics Institute, 2003.

\bibitem{HV}
K.~Hori and C.~Vafa.
\newblock Mirror symmetry.
\newblock {\em arXiv preprint hep-th/0002222}, 2000.

\bibitem{Iritani}
H.~Iritani.
\newblock Convergence of quantum cohomology by quantum {L}efschetz.
\newblock {\em J. Reine Angew. Math.}, 610:29--69, 2007.

\bibitem{Krawitz}
M.~Krawitz.
\newblock {FJRW} rings and {L}andau-{G}inzburg mirror symmetry.
\newblock {\em arXiv preprint arXiv:0906.0796}, 2009.

\bibitem{KS}
M.~Kreuzer and H.~Skarke.
\newblock On the classification of quasihomogeneous functions.
\newblock {\em Communications in Mathematical Physics}, 150:137--147, 1992.

\bibitem{Landman}
A.~Landman.
\newblock On the {P}icard-{L}efschetz transformation for algebraic manifolds
  acquiring general singularities.
\newblock {\em Trans. Amer. Math. Soc.}, 181:89--126, 1973.

\bibitem{LLY1}
B.~Lian, K.~Liu, and S.-T. Yau.
\newblock Mirror principle {I}.
\newblock {\em Asian J. of Math.}, 1:729--763, 1997.

\bibitem{LLY2}
B.~Lian, K.~Liu, and S.-T. Yau.
\newblock Mirror principle {II}.
\newblock {\em Asian J. of Math.}, 3:109--146, 1999.

\bibitem{LLY3}
B.~Lian, K.~Liu, and S.-T. Yau.
\newblock Mirror principle {III}.
\newblock {\em Asian J. of Math.}, 3:771--800, 1999.

\bibitem{LLY4}
B.~Lian, K.~Liu, and S.-T. Yau.
\newblock Mirror principle {IV}.
\newblock {\em Suv. Diff. Geom.}, VII:475--496, 2000.

\bibitem{McDuff}
D.~McDuff and D.~Salamon.
\newblock {\em Introduction to Symplectic Topology}.
\newblock Oxford Mathematical Monographs. Oxford University Press, 1999.

\bibitem{Morrison2}
D.~Morrison.
\newblock Mirror symmetry and rational curves on quintic threefolds: a guide
  for mathematicians.
\newblock {\em Jour. AMS}, 6:223--247, 1993.

\bibitem{Reid}
M.~Reid.
\newblock Decomposition of toric morphisms.
\newblock In M.~Artin and J.~Tate, editors, {\em Arithmetic and Geometry, vol.
  {II}}, volume~36 of {\em Progr. Math.}, pages 395--418. Birkh\"auser, 1983.

\bibitem{Schmid}
W.~Schmid.
\newblock Variation of {H}odge structure: the singularities of the period
  mapping.
\newblock {\em Inv. math.}, 22:211--319, 1973.

\bibitem{Wi93b}
E.~Witten.
\newblock Phases of {N} = 2 theories in two dimensions.
\newblock {\em Nuclear Phys. B}, 403(1):159--222, 1993.

\end{thebibliography}
\nocite{*}

\end{document}